\documentclass[preprint,12pt]{elsarticle}

\usepackage{amsmath, amsthm, amssymb}
\usepackage{epstopdf}
\usepackage{graphicx}
\usepackage{psfrag}

\usepackage{url}
\newcommand{{\Reals}}{\mathbb{R}}

\begin{document}

\begin{frontmatter}

\title{Fluid-structure interaction in blood flow capturing non-zero longitudinal structure displacement\tnoteref{label1}}
 \tnotetext[label1]{}

 \author[lab1]{Martina Buka\v{c}}
 \ead{martina@math.uh.edu}
% \cortext[cor1]{cor1}
 \author[lab1]{Sun\v{c}ica \v{C}ani\'{c}}
 \ead{canic@math.uh.edu}
% \fntext[label2]{label2}
% \cortext[cor1]{cor1}
 \author[lab1,lab2]{Roland Glowinski} %\corref{cor1}\fnref{label2}}
 \ead{roland@math.uh.edu}
% \cortext[cor1]{cor1}
 \author[lab3]{Josip Tamba\v{c}a} %\corref{cor1}\fnref{label2}}
 \ead{tambaca@math.hr}
 %\cortext[cor1]{cor1}
\author[lab1]{Annalisa Quaini} %\corref{cor1}\fnref{label2}}
 \ead{quaini@math.uh.edu}
% \cortext[cor1]{cor1}
 \address[lab1]{Department of Mathematics, University of Houston, 4800 Calhoun Rd, Houston, TX 77204, USA}
 \address[lab2]{Laboratoire Jacques-Louis Lions, Universit\'e P. et M. Curie, 4 Place Jussieu, 75005 Paris, France}
 \address[lab3]{Department of Mathematics, University of Zagreb, Bijeni\v{c}ka 30, 10000 Zagreb, Croatia}


%\author{}

%\address{}

\begin{abstract}
We present a new model and a novel loosely coupled partitioned numerical scheme 
modeling fluid-structure interaction (FSI) in blood flow allowing non-zero longitudinal displacement. 
Arterial walls are modeled by a {linearly viscoelastic, cylindrical Koiter shell model capturing both radial and longitudinal displacement}. 
Fluid flow is modeled by the Navier-Stokes equations for an incompressible, viscous fluid.
The two are fully coupled via kinematic and dynamic coupling conditions.
Our numerical scheme is based on 
a new modified Lie operator splitting that decouples the fluid and structure sub-problems in a
way that leads to a loosely coupled scheme which is { unconditionally} stable.
This was achieved by a clever use of the kinematic coupling condition at the fluid and structure sub-problems, 
leading to an implicit coupling between the fluid and structure velocities.
The proposed scheme is a modification of the recently introduced ``kinematically coupled scheme''
for which the newly proposed modified Lie splitting significantly increases the accuracy. 
The performance and accuracy of the scheme were studied on a couple of instructive examples
including a comparison with a monolithic scheme.
It was shown that 
the accuracy of our scheme was comparable to that of the monolithic scheme,
while our scheme retains all the main advantages 
of partitioned schemes, such as modularity, simple implementation, and low computational costs.
%Using energy-based arguments we showed that our scheme is unconditionally stable.
\end{abstract}


\end{frontmatter}

\section{Introduction}
We study fluid-structure interaction (FSI) between an incompressible viscous, Newtonian fluid, and a thin viscoelastic structure modeled
by the linearly viscoelastic cylindrical Koiter shell model.
The cylindrical viscoelastic Koiter shell model is derived to describe the mechanical properties of arterial walls, 
while the Navier-Stokes equations for an incompressible,
viscous, Newtonian fluid were employed to model the flow of blood in medium-to-large human arteries.
The two are coupled via the kinematic (no-slip) and dynamic (balance of contact forces) coupling conditions.
Motivated by recent results of {\sl in vivo} measurements of arterial wall motion~\cite{cinthio2006longitudinal,cinthio2005evaluation,persson2003new,svedlund2011longitudinal}, which 
indicate that both the radial and longitudinal
displacement, as well as viscoelasticity of arterial walls, are important in disease formation, we derived in this work the 
{viscoelastic cylindrical Koiter shell model} which captures both 
radial {and longitudinal displacement},
with the viscoelasticity of Kelvin-Voigt type.
The novel Koiter shell model is then coupled to the Navier-Stokes equations, and the coupled FSI problem is solved numerically.
In this manuscript we devise a stable, loosely coupled scheme to numerically solve the fully coupled FSI problem.
The scheme is based on a novel modified Lie's time-splitting, and on an implicit use of the kinematic coupling condition, as in~\cite{guidoboni2009stable},
which provides stability of the scheme without the need for sub-iterations between the 
fluid and structure sub-solvers. Stability of the scheme was proved in \cite{stability} on the same, simplified benchmark problem as in \cite{causin2005added}.
%We provide energy estimates which show that the scheme is  unconditionally stable. 
We provide numerical results which show that the scheme is first-order accurate in time. 
We compare our results with the monolithic scheme of Badia, Quaini, and Quarteroni~\cite{quaini2009algorithms,badia2008splitting} 
showing excellent agreement and comparable accuracy.

In hemodynamics, the coupling between fluid and structure is highly nonlinear due to the fact that
the fluid and structure densities are roughly the same,  
making the inertia of the fluid and structure roughly equal.
In this regime, classical partitioned loosely coupled (or explicit) numerical schemes, which are based on the fluid and structure sub-solvers,
have been shown to be intrinsically unstable~\cite{causin2005added}
due to the miss-match between the discrete energy dictated by the numerical scheme, and the continuous energy of the coupled problem.
This has been associated with the explicit role of  the ``added mass effect'', introduced and studied in~\cite{causin2005added}.
To rectify this problem,
the fluid and structure sub-solvers need to be iterated until the energy balance at the discrete level 
approximates well the energy of the continuous coupled problem.
The resulting strongly coupled partitioned scheme, however, gives rise to extremely high computational costs.

To get around these difficulties, several different loosely coupled algorithms have been proposed 
that modify the classical strategy in coupling the fluid and structure sub-solvers. 
The method proposed in~\cite{badia2008fluid} uses a simple membrane model for the structure that can be easily embedded into the fluid equations and appears as a generalized Robin boundary condition. 
In this way the original problem reduces to a sequence of fluid problems with a generalized Robin boundary condition that can be solved using only the fluid solver.
A similar approach was proposed in~\cite{nobile2008effective}, where the fluid and structure are split in the classical way, but the fluid and structure sub-problems were linked via novel
transmission (coupling) conditions that improve the convergence rate. 
Namely, a linear combination of the dynamic and kinematic interface conditions was used to artificially redistribute the fluid stress on the interface,
thereby avoiding the difficulty associated with the added mass effect.  

A different stabilization of the loosely coupled (explicit) schemes was proposed in~\cite{burman2009stabilization} which is based on {Nitsche's method}~\cite{hansbo2005nitsche} 
with a time penalty term giving $L^2$-control on the fluid force variations at the interface. 
We further mention the scheme proposed in~\cite{badia2009robin}, where Robin-Robin type preconditioner is combined with Krylov iterations for the solution of the interface system.

For completeness, we also mention several semi-implicit schemes. 
The schemes proposed in~\cite{fernandez2006projection,astorino2009added,astorino2009robin} separate the computation of fluid velocity from the coupled pressure-structure velocity system, 
thereby reducing the computational costs. Similar schemes,  derived from algebraic splitting, were proposed in~\cite{badia2008splitting,quaini2007semi}.
We also mention~\cite{murea2009fast} where an optimization problem is solved at each time-step to achieve continuity of stresses and continuity of velocity at the interface. 

In our work we deal with the problems associated with the added mass effect by: (1) employing the kinematic coupling condition implicitly 
in all the sub-steps of the splitting, as in the kinematically coupled scheme first introduced in~\cite{guidoboni2009stable};
(2) treating the fluid sub-problem together with the viscous part of the structure equations so that the structure inertia
appears in the fluid sub-problem (made possible by the kinematic coupling condition), giving rise to the energy estimates 
that mimic those in the continuous problem. In this step,
a portion of the fluid stress and the viscous part of the structure
equations are coupled weakly, and implicitly,
thereby adding dissipative effects to the fluid solver and contributing to the overall stability of the scheme
(although the scheme is stable even if viscoelasticity of the structure is neglected).
The modification of the Lie splitting introduced in this manuscript
uses the remaining portion of the normal fluid stress (the pressure) to explicitly load the structure 
in the elastodynamics equations, significantly increasing the accuracy of our scheme
when compared with the classical kinematically coupled scheme~\cite{guidoboni2009stable},
and making it comparable to that of the monolithic scheme presented in~\cite{quaini2009algorithms,badia2008splitting}. 
\begin{figure}[ht]
 \centering{
 \includegraphics[scale=0.6]{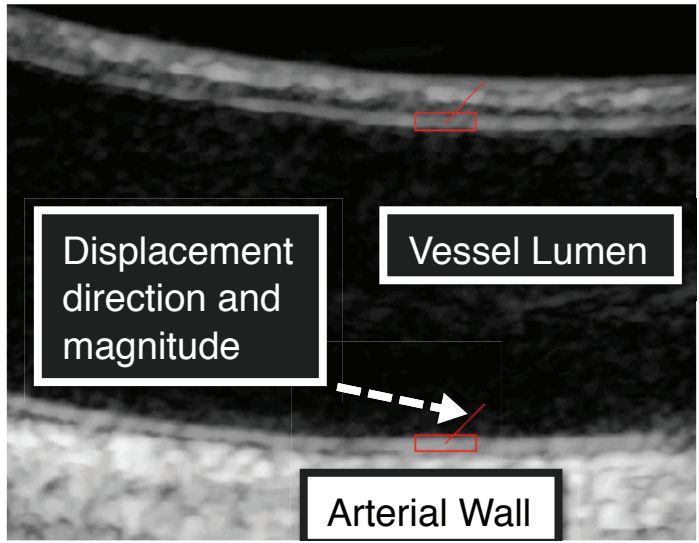}
 }\label{long_displ}
 \caption{Longitudinal displacement in a carotid artery measured using {\sl in vivo} ultrasound speckle tracking method. 
The thin red line located at the intimal layer of the arterial wall
 shows the direction and magnitude of the displacement vector, showing equal magnitude
in longitudinal and radial components of the displacement \cite{LyonLongitudinalMovie}.}
 \end{figure}

To deal with the motion of the fluid domain, we implemented an Arbitrary Lagrangian-Eulerian (ALE) approach.
In addition to the ALE method~\cite{hughes1981lagrangian,donea1983arbitrary,heil2004efficient,leuprecht2002numerical,quaini2007semi,quarteroni2000computational,le2001fluid,Bazilevs1,Bazilevs2}, the Immersed Boundary Method~\cite{fauci2006biofluidmechanics,fogelson2004platelet,peskin1977numerical,lim2004simulations,miller2005computational,peskin1989three}
has been very popular in problems with moving domains, especially when the structure is completely immersed in the fluid domain.
We also mention the Fictitious Domain Method combined with the mortar element method or ALE method~\cite{van2004combined,baaijens2001fictitious},
the Lattice Boltzmann method~\cite{fang2002lattice,feng2004immersed,krafczyk1998analysis,krafczyk2001two}, the Coupled Momentum Method~\cite{figueroa2006coupled}, and the Level Set Method~\cite{cottet2008eulerian}.
 
Even though other viscoelastic models have been proposed in literature to study FSI between blood and arterial walls, see e.g., \cite{Taylor,Pontrelli,canic2006modeling,Annals_viscoelastic},
to the best of our knowledge, the present manuscript is the first in which both radial and longitudinal displacement of a thin, 
viscoelastic structure are captured.
The main motivation for this work comes from recent developments in ultrasound speckle tracking techniques
(see \cite{cobbold2007foundations}, Chapter 8), which enabled {\sl in vivo} measurements of both the longitudinal and diameter vessel wall changes
over the cardiac cycle, indicating that longitudinal wall displacements can be comparable to the radial displacements, and should be 
included when studying tissue movement~\cite{warriner2008viscoelastic,cinthio2005evaluation,cinthio2006longitudinal}.
See Figure~\ref{long_displ}.
This is particularly pronounced under adrenaline 
conditions during which the longitudinal displacement of the intima-media 
complex increases by 200\%, and becomes twice the magnitude of radial displacement~\cite{ahlgren2009effects}.

The structure model capturing both radial and longitudinal displacement is presented next.

\section{The cylindrical linearly viscoelastic Koiter shell model}\label{sec2}

In this section we present the main steps in the derivation of the cylindrical linearly viscoelastic Koiter shell model that includes both
longitudinal and radial components of the displacement.

\begin{figure}[ht]
 \centering{
 \includegraphics[scale=0.5]{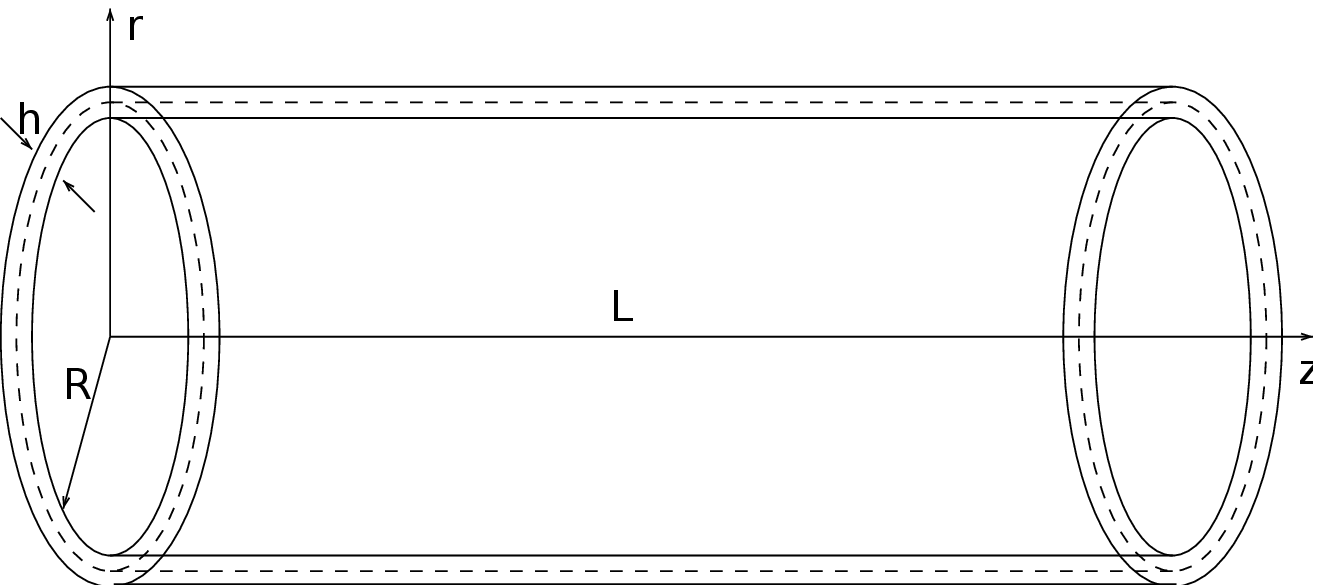}
 \includegraphics[scale=0.65]{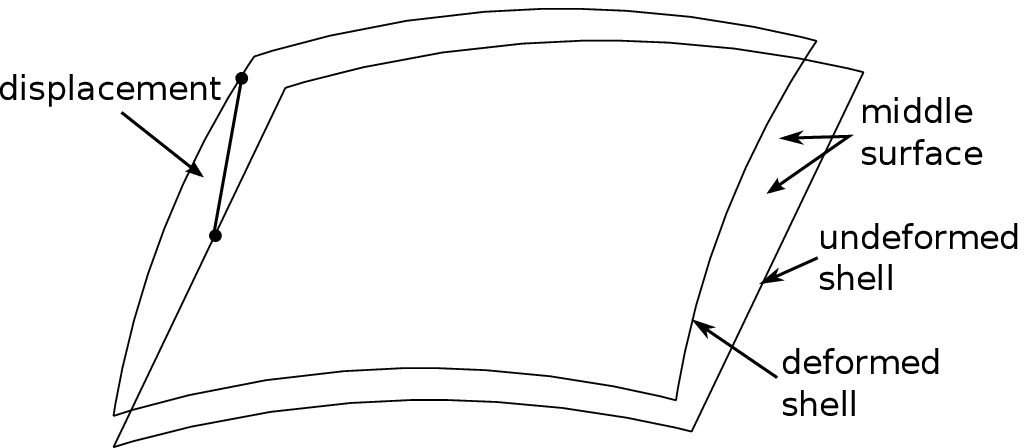}
 }\label{F1}
 \caption{Left: Cylindrical shell in reference configuration with middle surface radius $R$ and shell thickness $h$. Right: Deformed shell.}
 \end{figure}

Consider a clamped cylindrical shell of thickness $h$, length $L$, and reference radius of the middle surface equal to $R$.
See Figure~\ref{F1}. 
This reference configuration of the thin cylindrical shell will be denoted by  
\begin{equation}
\Gamma = \{ x=(R \cos \theta, R \sin \theta, z) \in \mathbb{R}^3 \; : \; \theta \in (0, 2 \pi), z \in (0,L)\}.
\end{equation} 
Displacement of the shell corresponds to the displacement of the shell's middle surface. 
Following the basic assumptions under which the Koiter shell model is valid~\cite{ciarlet1996asymptoticIII}, 
we assume that the walls of the cylinder are thin with respect to its radius (i.e., $h/R \ll 1$),  homogeneous, and deform linearly
(i.e., the displacement and the displacement gradient are both small). 
We will be assuming that the load exerted onto the shell is axially symmetric, leading to the axially symmetric displacements, 
so that the displacement in the $\theta-$direction is zero, and nothing in the problem depends on $\theta$.  
Thus, displacement $\boldsymbol\eta$ will have two components, the axial component $\eta_z$ and the radial component $\eta_r$.

To introduce the viscous effects into the linearly elastic Koiter shell model, we must assume that displacement is a function of both position and time. 
Thus $\boldsymbol\eta(z,t) = (\eta_z(z,t), \eta_r(z,t))$.
The derivatives with respect to the spatial variable will be denoted by $\boldsymbol\eta'$, and with respect to the temporal variable by $\dot{\boldsymbol\eta}$.

To define the Koiter shell model we need to define the geometry of deformation, and the physics of the problem which will be described by the
dynamic equilibrium equations (the Newton's second law of motion).

{\bf Geometry of Deformation.} The axially symmetric configuration of a Koiter shell can suffer stretching of the middle surface, and flexure (bending). 
The stretching of the middle surface is measured by the change of metric tensor, while flexure is measured by the change of curvature tensor.
The change of metric and the change of curvature tensors for a cylindrical shell are given, respectively by~\cite{ciarlet2000mathematical}
\begin{equation}
\boldsymbol{\gamma} (\boldsymbol\eta)=\left[ \begin{array}{cc} \eta_z' & 0 \\ 0 & R \eta_r \end{array} \right], \quad \boldsymbol\varrho (\boldsymbol\eta)=\left[ \begin{array}{cc} -\eta_r'' & 0 \\ 0 &  \eta_r \end{array} \right]. 
\end{equation}

{\bf The dynamic equilibrium equations.} The total energy of the linearly elastic Koiter shell is given by the sum of contributions due to 
stretching and flexure. The corresponding weak formulation will thus account for the 
internal (stretching) force and bending moment. 
For a linearly viscoelastic Koiter shell each of the contributions will also have an additional component to account for a viscoelastic response of the Koiter shell.
The elastic response will be modeled via the elasticity tensor $\mathbf{\mathcal{A}}$,
while the viscous response will be modeled by the viscosity tensor $\mathbf{\mathcal{B}}$.
Viscoelasticity will be modeled by the Kelvin-Voigt model 
in which the total stress is linearly proportional to strain and to the time-derivative of strain. 

More precisely,
introduce the {\bf elasticity} tensor $\mathbf{\mathcal{A}}$,~\cite{ciarlet2000mathematical},  to be defined as:
\begin{equation} 
 \mathbf{\mathcal{A}} \textbf{E}=\frac{2 E \sigma}{1-\sigma^2}(\textbf{A}^c \cdot \textbf{E})\textbf{A}^c+ \frac{2 E}{1+\sigma} \textbf{A}^c \textbf{E} \textbf{A}^c, \quad \textbf{E} \in \textrm{Sym}(\mathbb{R}^2),
\end{equation}
%or, in terms of the Lam\'{e} constants $\lambda$ and $\mu$:
%\begin{equation} 
 %\mathbf{\mathcal{A}} \textbf{E}=\frac{4\lambda\mu}{\lambda+2\mu}(\textbf{A}^c \cdot \textbf{E})\textbf{A}^c
 %+ 4\mu \textbf{A}^c \textbf{E} \textbf{A}^c, \quad \textbf{E} \in \textrm{Sym}(\mathbb{R}^2),
%\end{equation}
where $\textbf{A}_c=\left[ \begin{array}{cc} 1 & 0 \\ 0 & R^2 \end{array} \right]$ is the first fundamental form of the middle surface, $\textbf{A}^c=\left[ \begin{array}{cc} 1 & 0 \\ 0 &  \frac{1}{R^2} \end{array} \right] \textrm{is its inverse}$,
and $\cdot$ denotes the scalar product
$$A \cdot B := Tr(AB^{\tau}), \quad A,B \in M_2(\mathbb{R}).$$
Here $E$ is the \textit{Young's modulus} and $\sigma$ is the \textit{Poisson's ratio}.

Similarly, define the {\bf viscosity} tensor $\mathbf{\mathcal{B}}$ by:
\begin{equation}\label{vis}
\mathbf{\mathcal{B}} \textbf{E}=\frac{2 E_v \sigma_v}{1-\sigma_v^2}(\textbf{A}^c \cdot \textbf{E})\textbf{A}^c+ \frac{2 E_v}{1+\sigma_v} \textbf{A}^c \textbf{E} \textbf{A}^c, \quad \textbf{E} \in \textrm{Sym}(\mathbb{R}^2).
\end{equation}
Here $E_v$ and $\sigma_v$ correspond to the viscous counterparts of the Young's modulus $E$ and Poisson's ratio $\sigma$. 
Then, for a linearly viscoelastic Koiter shell model we define the internal (stretching) force
\begin{equation}
N := \frac{h}{2} \mathbf{\mathcal{A}} \boldsymbol\gamma(\boldsymbol\eta)+\frac{h}{2} \mathbf{\mathcal{B}} \boldsymbol\gamma(\dot{\boldsymbol\eta}),
\end{equation}
and bending moment
\begin{equation}
M := \frac{h^3}{24} \mathbf{\mathcal{A}} \boldsymbol\varrho(\boldsymbol\eta)+\frac{h^3}{24} \mathbf{\mathcal{B}} \boldsymbol\varrho(\dot{\boldsymbol\eta}).
\end{equation}
The weak formulation of the linearly viscoelastic Koiter shell is then given by the following: for each $t>0$ find $\boldsymbol\eta(\cdot,t) \in V_c$ such that $\forall \boldsymbol\xi \in V_c$
\begin{equation}\label{weakKoiter}
\begin{array}{c}
\displaystyle{\frac{h}{2} \int_0^L (\mathbf{\mathcal{A}} \boldsymbol\gamma(\boldsymbol\eta) + \mathbf{\mathcal{B}} \boldsymbol\gamma(\dot{\boldsymbol\eta})) \cdot \boldsymbol\gamma(\boldsymbol\xi) R dz + \frac{h^3}{24} \int_0^L (\mathbf{\mathcal{A}} \boldsymbol\varrho(\boldsymbol\eta) + \mathbf{\mathcal{B}} \boldsymbol\varrho(\dot{\boldsymbol\eta})) \cdot \boldsymbol\varrho(\boldsymbol\xi) R dz}\\
\\
+ \displaystyle{\rho_s h \int_0^L \frac{\partial^2 \boldsymbol\eta}{\partial t^2} \cdot \boldsymbol\xi R dz= \int_0^L \mathbf{f} \cdot \boldsymbol\xi R dz,}
\end{array}
\end{equation}
where $\rho_s$ denotes the volume shell density and 
\begin{eqnarray*}
 V_c &=& H_0^1(0,L) \times H_0^2(0,L) \\
        &=& \{\boldsymbol\xi = (\xi_z, \xi_r) \in H^1(0,L) \; : \; \boldsymbol\xi(0)=\boldsymbol\xi(L)=0, \; \xi_r'(0)=\xi_r'(L)=0\}.
\end{eqnarray*}
The first term on the left-hand side of~\eqref{weakKoiter} multiplying $h/2$ captures the membrane effects, while the second
term on the left hand-side multiplying $h^3/24$ captures the flexural effects of the Koiter shell.

Components of the forcing term $\mathbf{f} = (f_z,f_r)^T$ are the surface densities in the reference configuration of the axial and radial force.
The corresponding dynamic equilibrium equations in differential form can be written as follows:
\begin{eqnarray}
\rho_s h \frac{\partial^2 \eta_z}{\partial t^2}-C_2 \frac{\partial \eta_r}{\partial z}-C_3 \frac{\partial^2 \eta_z}{\partial z^2}-D_2 \frac{\partial ^2 \eta_r}{\partial t \partial z} -D_3 \frac{\partial ^3 \eta_z}{\partial t \partial z^2} = f_z  \label{structure1}  \\
\rho_s h \frac{\partial^2 \eta_r}{\partial t^2}+C_0 \eta_r -C_1 \frac{\partial^2 \eta_r}{\partial z^2}  +C_2 \frac{\partial \eta_z}{\partial z} + C_4 \frac{\partial^4 \eta_r}{\partial z^4}+D_0 \frac{\partial \eta_r}{\partial t}-D_1 \frac{\partial^3 \eta_r}{\partial t \partial z^2} \nonumber  \\
\quad +D_2 \frac{\partial ^2 \eta_z}{\partial t \partial z}+D_4 \frac{\partial ^5 \eta_r}{\partial t \partial z^4} = f_r, \label{structure2}
\end{eqnarray}
where 
\begin{equation}
\small{
 \begin{array}{rlrlrl}
 C_0 &= \displaystyle{\frac{h E}{R^2(1-\sigma^2)}(1+\frac{h^2}{12 R^2}),}  & C_1 
 &\displaystyle{= \frac{h^3}{6} \frac{E \sigma}{R^2 (1-\sigma^2)}, } & C_2 
 &\displaystyle{=\frac{h}{R}\frac{E \sigma}{1-\sigma^2}, } \\
C_3&\displaystyle{=\frac{h E}{1-\sigma^2},} & C_4 &\displaystyle{=\frac{h^3}{12}\frac{E}{1-\sigma^2},}\\
\\
D_0 &= \displaystyle{\frac{h}{R^2} C_v(1+\frac{h^2}{12 R^2}),}  & D_1 &=\displaystyle{ \frac{h^3}{6} \frac{D_v}{R^2}, } &
 D_2 &= \displaystyle{\frac{h D_v}{R}, }\\
 D_3 &= h C_v, & D_4 &=\displaystyle{ \frac{h^3}{12}C_v,}
\label{coeff}
\end{array}}
\end{equation}
and
$$C_v:= \frac{E_v}{1-\sigma_v^2}, \quad D_v:= \frac{E_v \sigma_v}{1-\sigma_v^2}.$$
The boundary conditions for a clamped Koiter shell problem are given by
\begin{equation}\label{boundaryKoiter}
\boldsymbol\eta(0,t)=\boldsymbol\eta(L,t)=0,\quad \frac{\partial\eta_r}{\partial z}(0,t)=\frac{\partial\eta_r}{\partial z}(L,t)=0.
\end{equation}
For a mathematical justification of the Koiter shell model please see \cite{ciarlet1996asymptoticIII,Xiao,Li,Li2,Li3}.
  
\section{The fluid-structure interaction problem}\label{sec3}

\begin{figure}[ht]
\centering{
\includegraphics[scale=0.7]{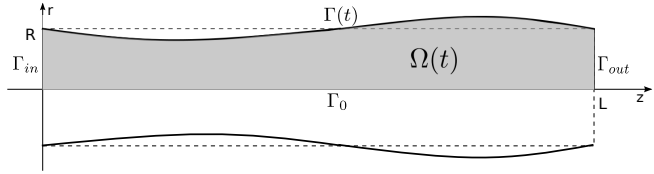}
}
\caption{Deformed domain $\Omega(t).$}
\label{fig:domain}
\end{figure}

We consider the flow of an incompressible, viscous fluid in a two-dimensional channel of reference length $L$, and reference width $2R$,
see Figure~\ref{fig:domain}.
The lateral boundary of the channel is bounded by a thin, deformable wall, modeled by
the linearly viscoelastic Koiter shell model, described in the previous section.
We are interested in simulating a pressure-driven flow through the deformable 2D channel with a two-way coupling 
between the fluid and structure.
Without loss of generality, we consider only the upper half of the fluid domain supplemented by a symmetry condition at the axis of symmetry.
Thus, the reference domain in our problem is given by
\begin{equation*}
\Omega_0 := \{(z,r) | 0<z<L, 0<r<R \}.
\end{equation*}
Here $z$ and $r$ denote the horizontal and vertical Cartesian coordinates, respectively.
See Figure~\ref{fig:domain}.

\vskip 0.1in
\noindent
{\bf Remark 1.} It is worth mentioning here that while the fluid flow will be modeled in 2D, the thin structure equations, 
described in the previous section, are given in terms of cylindrical coordinates, assuming axial symmetry. It is standard practice in 
2D fluid-structure interaction studies to use thin structure equations that are 
derived assuming cylindrical geometry. This is because
cylindrical structure models
account for the circumferential stress that ``keeps" the top and bottom boundary of the structure ``coupled together"
when they are loaded by the stresses exerted by the fluid, thereby giving rise to physiologically reasonable solutions.
\vskip 0.1in

The mathematical model for the corresponding fluid-structure interaction problem can be defined as follows.
The fluid domain, which depends on time, is not known {\sl a priori}. 
%The lateral boundary of the fluid domain is given by the motion of the linearly viscoelastic Koiter shell.
The location of the lateral boundary, defined in Lagrangian framework, is given by
$\Gamma(t) = \{(\hat{z}+\eta_z(\hat{z},t),R+\eta_r(\hat{z},t))\ | \ \hat{z}\in(0,L)\}$ for $t \in (0,T)$.
Throughout the rest of the manuscript we will be denoting the Lagrangian coordinates by $\hat{\mathbf{x}}=(\hat{z},\hat{r})$. The displacement of the boundary will always be given in Lagrangian framework. However, we will omit the hat notation on $\boldsymbol \eta$ for simplicity. 
%As mentioned earlier, we will be assuming that the shell is fixed at the end-points so that $\eta_z(L,t)=\eta_r(L,t)= 0$.

We will be assuming that for each $t\in (0,T)$ the boundary of the fluid domain is Lipschitz continuous (see Figure \ref{fig:domain}), and that it's lateral boundary,
in Eulerian framework, can be described by a Lipschitz continuous function
\begin{equation*}
g(\cdot\ ;t) : (0,L)\to \mathbb{R}, \quad g(\cdot\ ;t) : z \mapsto g(z;t),\quad {\rm for\ each} \quad t\in(0,T)
\end{equation*}
so that, in Eulerian framework,
\begin{equation*}
\Gamma(t) = \{(z,g(z;t)), z \in (0,L)\}, \quad {\rm for}\quad t\in(0,T).
\end{equation*}
The fluid domain is given by
\begin{equation}
 \Omega(t) = \{(z,r) \in \mathbb{R}^2; \; 0<z<L, \; 0< r< g(z;t)  \},\quad {\rm for}\quad t\in(0,T).
\end{equation}
The inlet boundary will be denoted by $\Gamma_{\rm in}$, the outlet boundary by $\Gamma_{\rm out}$, 
the symmetry (bottom) boundary for which $r = 0$ by $\Gamma_0$, so that 
\begin{equation*}
\partial\Omega(t) = \Gamma_{\rm in}\cup \Gamma(t) \cup \Gamma_{\rm out}\cup \Gamma_0.
\end{equation*}
%See Figure~\ref{fig:domain}.

The flow of a viscous, incompressible, Newtonian fluid is governed by the Navier-Stokes equations
\begin{eqnarray}\label{NS1}
 \rho_f \bigg( \frac{\partial \mathbf{u}}{\partial t}+ \mathbf{u} \cdot \nabla \mathbf{u} \bigg) & =& \nabla \cdot \boldsymbol\sigma \quad \;  \textrm{in}\; \Omega(t)\; \textrm{for}\; t \; \in(0,T), \\
 \label{NS2}
\nabla \cdot \mathbf{u} &=& 0 \quad \quad \; \quad \textrm{in}\; \Omega(t)\; \textrm{for}\; t \; \in(0,T),
\end{eqnarray}
where $\mathbf{u}=(u_z,u_r)$ is the fluid velocity, $p$ is the fluid pressure, $\rho_f$ is the fluid density, and $\boldsymbol\sigma$ is the fluid stress tensor. 
For a Newtonian fluid the stress tensor is given by $\boldsymbol\sigma = -p \mathbf{I} + 2 \mu \mathbf{D}(\mathbf{u}),$ where
$\mu$ is the fluid viscosity and  $\mathbf{D}(\mathbf{u}) = (\nabla \mathbf{u}+(\nabla \mathbf{u})^{\tau})/2$ is the rate-of-strain tensor.

At the inlet and outlet boundary we prescribe the normal stress:
\begin{eqnarray}
\mathbf{\boldsymbol\sigma n}|_{in}(0,r,t) &=& -p_{in}(t) \mathbf{n}|_{in} \quad \textrm{on} \; (0,R) \times (0,T), \label{inlet} \\
\mathbf{\boldsymbol\sigma n}|_{out}(L,r,t) &=& -p_{out}(t) \mathbf{n}|_{out} \quad \textrm{on} \; (0,R) \times (0,T),\label{outlet}
\end{eqnarray}
where $\mathbf{n}_{in}/\mathbf{n}_{out}$ are the outward normals to the inlet/outlet boundaries, respectively.
These boundary conditions are common in blood flow modeling~\cite{badia2008fluid,miller2005computational,nobile2001numerical}.

At the bottom boundary $r=0$ we impose the  symmetry conditions:
\begin{equation}\label{symmetry_condition}
 \frac{\partial u_z}{\partial r}(z,0,t) = 0, \quad u_r(z,0,t) = 0 \quad \textrm{on} \; (0,L) \times (0,T).
\end{equation}

The upper boundary $\Gamma(t)$ represents the deformable channel wall, whose dynamics is modeled by~\eqref{structure1}-\eqref{structure2}.
The structure equations are supplemented with the following boundary conditions
\begin{equation}\label{homostructure}
 \boldsymbol \eta(0,t) = \boldsymbol \eta(L,t) = \partial_z \eta_r(0,t)=  \partial_z \eta_r(L,t) = 0 \quad \textrm{on} \; (0,T).
\end{equation}

\if 1  0
\vskip 0.1in
\noindent
{\bf Remark 3.}
Although not physiologically optimal,
the homogeneous boundary conditions~\eqref{homostructure} ``contaminate'' the flow only locally,
near the inlet and outlet boundaries, as it was shown in~\cite{canic2003effective} that the boundary layer that is formed
due to the homogeneous boundary conditions decays exponentially fast away from 
the boundary. 
\vskip 0.1in
\fi

Initially, the fluid and the structure are assumed to be at rest, with zero displacement from the reference configuration
\begin{equation}\label{initial}
 \mathbf{u}=0, \quad \boldsymbol\eta = 0, \quad \frac{\partial \boldsymbol\eta}{\partial t}=0.
 \end{equation}

The fluid and structure are coupled via the kinematic and dynamic boundary conditions~\cite{canic2005self}:
\begin{itemize}
 \item \textbf{Kinematic coupling condition} 
describes continuity of velocity
\begin{equation}\label{kinematic}
 \mathbf{u}(\hat{z}+\eta_z(\hat{z},t),R+\eta_r(\hat{z},t),t)=\frac{\partial \boldsymbol\eta}{\partial t}(\hat{z},t) \quad \; \textrm{on} \; (0,L)\times (0,T).
\end{equation} 
\item \textbf{Dynamic coupling condition} 
describes balance of contact forces, namely, it says that the contact force exerted by the fluid is equal but of opposite sign to the contact force exerted by the structure to the fluid:
\begin{eqnarray}
 f_z &=& - J\ \boldsymbol {\widehat{\sigma \mathbf{n}}}|_{\Gamma(t)}  \cdot \mathbf{e_z}  \; \textrm{on} \; (0,L)\times(0,T),   \label{dynamic1}\\
 f_r &=& - J\  \boldsymbol {\widehat{\sigma \mathbf{n}}}|_{\Gamma(t)} \cdot \mathbf{e_r}    \; \textrm{on} \; (0,L)\times (0,T),
\label{dynamic2}
\end{eqnarray}
where
\begin{equation}\label{Jacobian}
J = \sqrt{\bigg(1+ \frac{\partial \eta_z}{\partial z} \bigg)^2+\bigg(\frac{\partial \eta_r}{\partial z}\bigg)^2}
\end{equation}
denotes the Jacobian of transformation from the Eulerian to the Lagrangian framework, and $\boldsymbol {\widehat{\sigma \mathbf{n}}}$ denotes the normal fluid stress on the reference domain $\hat{\Omega} = (0,L)\times(0,R)$. Here $\mathbf{e_z}=(1,0)$ and $ \mathbf{e_r}=(0,1)$ are the standard unit basis vectors, and $\boldsymbol n$ is the outward normal to the deformed domain.
\end{itemize}

\subsection{The energy of the coupled FSI problem}
To formally derive the energy of the coupled FSI problem we multiply the structure equations by the structure velocity,
the balance of momentum in the fluid equations by the fluid velocity, integrate by parts over the respective domains
using the incompressibility condition, 
and add the two equations together.
The dynamic and kinematic coupling conditions are then used to couple the fluid and structure sub-problems.
The resulting equation represents the total energy of the problem.

We start by first considering the Koiter shell model for the structure. We recall the weak formulation
of the clamped Koiter shell given by \eqref{weakKoiter}.
In \eqref{weakKoiter} we replace the test function $\boldsymbol\xi$ by the structure velocity $\frac{\partial \boldsymbol \eta}{\partial t}$
and integrate by parts over $(0,L)$ to obtain the following energy equality of the clamped Koiter shell:
 \begin{eqnarray}
& \frac{d}{dt} \bigg\{
 \frac{\rho_s h}{2} \bigg|\bigg|\frac{\partial \eta_z}{\partial t} \bigg|\bigg|^2_{L^2(0,L)} + \frac{\rho_s h}{2}  \bigg|\bigg|\frac{\partial \eta_r}{\partial t} \bigg|\bigg|^2_{L^2(0,L)}
 \nonumber\\
& +\frac{h}{2} \left[
  \frac{E}{1+\sigma}\bigg|\bigg|\frac{\eta_r}{R}\bigg|\bigg|^2_{L^2(0,L)} +
  \frac{E}{1+\sigma}  \bigg|\bigg|\frac{\partial \eta_z}{\partial z} \bigg|\bigg|^2_{L^2(0,L)}
  +\frac{E \sigma}{1-\sigma^2}  \bigg|\bigg|\frac{\partial \eta_z}{\partial z}+ \frac{\eta_r}{R} \bigg|\bigg|^2_{L^2(0,L)}
 \right]
 \nonumber
 \\
& +\frac{h^3}{24} \left[
 \frac{E}{1+\sigma}\bigg|\bigg|\frac{\eta_r}{R^2}\bigg|\bigg|^2_{L^2(0,L)} 
+  \frac{E}{1+\sigma}  \bigg|\bigg|\frac{\partial^2 \eta_r}{\partial z^2} \bigg|\bigg|^2_{L^2(0,L)}
 +\frac{E \sigma}{1-\sigma^2}  \bigg|\bigg|-\frac{\partial^2 \eta_r}{\partial z^2}+ \frac{\eta_r}{R^2} \bigg|\bigg|^2_{L^2(0,L)} 
 \right]
\bigg\}
\nonumber
\\
& +\frac{h}{2} \left[
  \frac{E_v}{1+\sigma_v}\bigg|\bigg|\frac{\partial\eta_r}{R\partial t}\bigg|\bigg|^2_{L^2(0,L)} +
  \frac{E_v}{1+\sigma_v}  \bigg|\bigg|\frac{\partial^2 \eta_z}{\partial z\partial t} \bigg|\bigg|^2_{L^2(0,L)}
  +\frac{E_v \sigma_v}{1-\sigma_v^2}  \bigg|\bigg|\frac{\partial^2 \eta_z}{\partial z\partial t}+ \frac{\partial\eta_r}{R\partial t} \bigg|\bigg|^2_{L^2(0,L)}
 \right]
 \nonumber
 \\
& +\frac{h^3}{24} \left[
 \frac{E_v}{1+\sigma_v}\bigg|\bigg|\frac{\partial\eta_r}{R^2\partial t}\bigg|\bigg|^2_{L^2(0,L)} 
+  \frac{E_v}{1+\sigma_v}  \bigg|\bigg|\frac{\partial^3 \eta_r}{\partial z^2\partial t} \bigg|\bigg|^2_{L^2(0,L)}
 +\frac{E_v \sigma_v}{1-\sigma_v^2}  \bigg|\bigg|-\frac{\partial^3 \eta_r}{\partial z^2\partial t}+ \frac{\partial\eta_r}{R^2\partial t} \bigg|\bigg|^2_{L^2(0,L)} 
 \right]
 \nonumber
 \\
 &= \int_0^L \mathbf{f}\cdot \frac{\partial \boldsymbol \eta}{\partial t}\  d\hat{z}.
 \label{KoiterEnergy}
\end{eqnarray}
The first two terms under the time-derivative correspond to the kinetic energy of the Koiter shell. 
%which we denote by $\cal{{E}}_{\rm s}^{\rm kin}$. 
The terms in the second and third row correspond to the elastic energy of the Koiter shell
% which we denote by$\cal{{E}}_{\rm s}^{\rm el}$
(the terms multiplying $h$ are the membrane energy, while the terms multiplying $h^3$ correspond to the flexural (bending) energy).
The terms in the fourth and fifth row correspond to the viscous energy of the viscoelastic Koiter shell, 
%which we denote by $\cal{{E}}_{\rm s}^{\rm vis}$.
while the last term corresponds to the work done by the external loading, which comes form the fluid stress via the dynamic coupling condition.

To deal with the fluid sub-problem we multiply the momentum equation in the Navier-Stokes equations by $\boldsymbol u$ and integrate by parts over $\Omega(t)$,
using the incompressibility condition along the way. 
With the help of the following identities
\begin{eqnarray*}
\int_{\Omega(t)} \frac{\partial \boldsymbol u}{\partial t} \boldsymbol u d\boldsymbol x &=& \frac{1}{2} \frac{d}{dt}\int_{\Omega(t)} |\boldsymbol u|^2 d\boldsymbol x - \frac{1}{2} \int_{\partial \Omega(t)} |\boldsymbol u|^2 \boldsymbol u \cdot \boldsymbol n dS, \label{ident1} \\
\int_{\Omega(t)} (\boldsymbol u \cdot \nabla) \boldsymbol u \cdot \boldsymbol u d\boldsymbol x &=& \frac{1}{2} \int_{\partial \Omega(t)} |\boldsymbol u|^2 \boldsymbol u \cdot \boldsymbol n dS, \label{ident2}
\end{eqnarray*}
one obtains
\begin{eqnarray*}
 \frac{1}{2}\frac{d}{dt} \bigg\{  \rho_f ||\boldsymbol u||^2_{L^2(\Omega(t))} \bigg\} + 2 \mu ||\boldsymbol D(\boldsymbol u)||^2_{L^2(\Omega(t))}\\
- \int_0^R p_{in}(t) u_z|_{z=0} dr +   \int_0^R p_{out}(t) u_z|_{z=L} dr = \int_{\Gamma(t)} \boldsymbol\sigma \mathbf{n} \cdot \boldsymbol u\ dS
\end{eqnarray*}
The integral on the right-hand side can be written in Lagrangian coordinates as
\begin{equation}
\int_{\Gamma(t)} \boldsymbol\sigma \mathbf{n} \cdot \boldsymbol u\ dS
=\int_0^L \left[\boldsymbol\sigma \mathbf{n} \cdot \boldsymbol u\right]|_{(\hat{z}+\eta_z(\hat{z},t),R+\eta_r(\hat{z},t))} \ J \ d\hat{z}
\end{equation}
where $J$ is the Jacobian of transformation from the Eulerian to Lagrangian framework, given by \eqref{Jacobian}.
Now we use the kinematic and dynamic  lateral boundary conditions \eqref{kinematic}-\eqref{dynamic2} to obtain
\begin{equation}
\int_0^L \left[\boldsymbol\sigma \mathbf{n} \cdot \boldsymbol u\right]|_{(\hat{z}+\eta_z(\hat{z},t),R+\eta_r(\hat{z},t))} \ J \ d\hat{z}
= -  \int_0^L \mathbf{f}\cdot \frac{\partial \boldsymbol \eta}{\partial t}\  d\hat{z}.
\end{equation}
Thus, the fluid sub-problem coupled with the structure satisfies
\begin{eqnarray}
 &\frac{1}{2}\frac{d}{dt} \bigg\{  \rho_f ||\boldsymbol u||^2_{L^2(\Omega(t))} \bigg\} + 2 \mu ||\boldsymbol D(\boldsymbol u)||^2_{L^2(\Omega(t))}
 \nonumber\\
 \nonumber\\
&- \displaystyle{\int_0^R p_{in}(t) u_z|_{z=0} dr +   \int_0^R p_{out}(t) u_z|_{z=L} dr} = -  \int_0^L \mathbf{f}\cdot \frac{\partial \boldsymbol \eta}{\partial t}\  d\hat{z}.
\label{FluidEnergy}
\end{eqnarray}
By adding \eqref{KoiterEnergy} and \eqref{FluidEnergy}, the right-hand sides
of the two equations cancel out and one obtains the energy equality
for the FSI problem:
{\scriptsize{
 \begin{eqnarray}
& \frac{d}{dt} \Bigg\{ \underbrace{ \frac{\rho_f}{2} ||\boldsymbol u||^2_{L^2(\Omega(t))} 
+  \frac{\rho_s h}{2} \bigg|\bigg|\frac{\partial \eta_z}{\partial t} \bigg|\bigg|^2_{L^2(0,L)} + \frac{\rho_s h}{2}  \bigg|\bigg|\frac{\partial \eta_r}{\partial t} \bigg|\bigg|^2_{L^2(0,L)}}_{Kinetic \ Energy}
 \nonumber\\
& +\frac{h}{2} \underbrace{\left[
  \frac{E}{1+\sigma}\bigg|\bigg|\frac{\eta_r}{R}\bigg|\bigg|^2_{L^2(0,L)} +
  \frac{E}{1+\sigma}  \bigg|\bigg|\frac{\partial \eta_z}{\partial z} \bigg|\bigg|^2_{L^2(0,L)}
  +\frac{E \sigma}{1-\sigma^2}  \bigg|\bigg|\frac{\partial \eta_z}{\partial z}+ \frac{\eta_r}{R} \bigg|\bigg|^2_{L^2(0,L)}\right]}_{Structure \ Elastic \ Energy\ (Membrane\ Contribution)}
 \nonumber
\\
& +\frac{h^3}{24} \underbrace{\left[
 \frac{E}{1+\sigma}\bigg|\bigg|\frac{\eta_r}{R^2}\bigg|\bigg|^2_{L^2(0,L)} 
+  \frac{E}{1+\sigma}  \bigg|\bigg|\frac{\partial^2 \eta_r}{\partial z^2} \bigg|\bigg|^2_{L^2(0,L)}
 +\frac{E \sigma}{1-\sigma^2}  \bigg|\bigg|-\frac{\partial^2 \eta_r}{\partial z^2}+ \frac{\eta_r}{R^2} \bigg|\bigg|^2_{L^2(0,L)}\right]}_{Structure \ Elastic\ Energy\ (Flexural(Shell) \ Contribution)}
\Bigg\}
\nonumber
\end{eqnarray}
\begin{eqnarray}
& +\frac{h}{2} \underbrace{\left[
  \frac{E_v}{1+\sigma_v}\bigg|\bigg|\frac{\partial\eta_r}{R\partial t}\bigg|\bigg|^2_{L^2(0,L)} +
  \frac{E_v}{1+\sigma_v}  \bigg|\bigg|\frac{\partial^2 \eta_z}{\partial z\partial t} \bigg|\bigg|^2_{L^2(0,L)}
  +\frac{E_v \sigma_v}{1-\sigma_v^2}  \bigg|\bigg|\frac{\partial^2 \eta_z}{\partial z\partial t}
+ \frac{\partial\eta_r}{R\partial t} \bigg|\bigg|^2_{L^2(0,L)}\right]}_{Structure\ Viscous \ Energy\ (Membrane\ Contribution)}
 \nonumber
 \\
& +\frac{h^3}{24} \underbrace{\left[
 \frac{E_v}{1+\sigma_v}\bigg|\bigg|\frac{\partial\eta_r}{R^2\partial t}\bigg|\bigg|^2_{L^2(0,L)} 
+  \frac{E_v}{1+\sigma_v}  \bigg|\bigg|\frac{\partial^3 \eta_r}{\partial z^2\partial t} \bigg|\bigg|^2_{L^2(0,L)}
 +\frac{E_v \sigma_v}{1-\sigma_v^2}  \bigg|\bigg|-\frac{\partial^3 \eta_r}{\partial z^2\partial t}
+ \frac{\partial\eta_r}{R^2\partial t} \bigg|\bigg|^2_{L^2(0,L)}\right]}_{Structure\ Viscous \ Energy\ (Flexural\ (Shell) \ Contribution)}
 \nonumber
 \\
 & +  \underbrace{2 \mu ||\boldsymbol D(\boldsymbol u)||^2_{L^2(\Omega(t))}}_{Fluid \ Viscous \ Energy} = 
  \displaystyle{\int_0^R p_{in}(t) u_z|_{z=0} dr -   \int_0^R p_{out}(t) u_z|_{z=L} dr}
 \label{FSIEnergy}
\end{eqnarray}
}}
The coefficients $E_v,\sigma_v$ are defined after equation~\eqref{vis}. Therefore, we have shown that
if a solution to the coupled fluid-structure interaction problem~\eqref{structure1} - \eqref{dynamic2} exists, then it satisfies 
the  energy equality~\eqref{FSIEnergy}. This equality says that the rate of change of the kinetic energy of the fluid,  the kinetic energy of the structure,
and  the elastic energy of the structure, plus the viscous energy of the structure, plus the viscous energy of the fluid,
is equal to the work done by the inlet and outlet data.

\section{The numerical scheme}\label{sec4}
To solve the fluid-structure interaction problem~\eqref{structure1}-\eqref{dynamic2},
we propose an extension of a loosely coupled partitioned scheme, called the kinematically coupled scheme, first introduced in~\cite{guidoboni2009stable}.

The classical kinematically coupled scheme introduced in~\cite{guidoboni2009stable}
is based on a time-splitting approach known as the Lie splitting~\cite{glowinski2003finite}.
The viscoelastic structure is split into its elastic part and the viscous part.
The viscous (parabolic) part is treated together with the fluid, while the elastic (hyperbolic) part is treated separately.
The inclusion of the viscous part of the structure into the fluid solver as a boundary condition in the weak formulation
accounts for the fluid load onto the structure as it simultaneously
deals with the ``added mass effect''~\cite{causin2005added}.
This adds dissipative effects to the fluid solver that contribute to the stability of the scheme. This approach provides a desirable
discrete energy inequality making this scheme stable even when the density of the fluid is equal to the density of the structure,
which is the case in the blood flow application. 
The elastic part of the structure,  which is solved separately, communicates with the fluid only via the kinematic coupling condition
in the classical kinematically coupled scheme.
The fluid stress does not appear in this step, as it is used as a loading to the viscous part of the structure
in the weak formulation of the fluid sub-problem.

In this manuscript we will change this approach by additionally splitting the normal stress into a fraction that loads the
viscous part of the structure, and a fraction (pressure) that loads the elastic part of the structure.  
%The fraction of the stress that loads the viscoelastic part of the structure will be included in an implicit way,
%subtracted by the (explicit) contribution of a fraction of 
%the normal stress (pressure) which will then be use to load 
%the elastic part of the structure.
This splitting is done using a modification of the Lie splitting scheme 
in a way which significantly increasases accuracy.

Thus, in this manuscript, the kinematically coupled scheme is extended and improved to achieve the following two goals:
\begin{enumerate}
\item Capture both the radial and longitudinal displacement of the linearly viscoelastic Koiter shell for the underlying fluid-structure interaction problem.
%This is an extension of the scheme presented in \cite{guidoboni2009stable} where longitudinal displacement was neglected. 
\item Increase the accuracy of the kinematically coupled scheme by introducing 
a new splitting strategy based on a modified Lie's scheme.
%In the new splitting strategy the normal fluid stress is split so that a fraction of it is coupled with the viscous part of the structure weakly and in
%an implicit way,
%while the remaining part is coupled with the elastodynamics of the structure strongly and in an explicit way.
\end{enumerate}

This version of the kinematically coupled scheme retains all the advantages of the original scheme, which include:
\begin{itemize}
\item The scheme does not require sub-iterations between the fluid and structure
sub-solvers to achieve stability. 
%This is due to the fact that the problems associated with the added mass effect are resolved elegantly by adding
%the contribution of the structure kinetic energy to the fluid energy making the total discrete energy bounded in time.
\item The scheme is modular, allowing the use of one's favorite fluid or structure solvers independently. 
%The choice of solvers 
%can then be easily tailored to capture
%different physical phenomena, such as, viscous dissipation (parabolic) vs. wave propagation (hyperbolic). 
The solvers communicate 
through the initial conditions.
\item Except for the pressure, the fluid stress at the boundary does not need to be calculated explicitly. 
\end{itemize}

\subsection{The Lie scheme}
To apply the Lie splitting scheme the problem must first be written as a first-order system in time:
\begin{eqnarray}\label{LieProblem}
   \frac{\partial \phi}{\partial t} + A(\phi) &=& 0, \quad \textrm{in} \ (0,T), \\
\phi(0) &=& \phi_0, 
\end{eqnarray}
where A is an operator from a Hilbert space into itself. Operator A is then split, in a non-trivial decomposition, as
\begin{equation}
 A = \sum\limits_{i=1}^I A_i.
\end{equation}
The Lie scheme consists of the following. Let $\triangle t>0$ be a time discretization step. Denote $t^n=n\triangle t$ and let $\phi^n$
be an approximation of $\phi(t^n).$ Set $\phi^0=\phi_0.$ Then, for $n \geq 0$ compute $\phi^{n+1}$ by solving
\begin{eqnarray}
   \frac{\partial \phi_i}{\partial t} + A_i(\phi_i) &=& 0 \quad \textrm{in} \; (t^n, t^{n+1}), \\
\phi_i(t^n) &=& \phi^{n+(i-1)/I}, 
\end{eqnarray}
and then set $\phi^{n+i/I} = \phi_i(t^{n+1}),$ for $i=1, \dots. I.$

This method is first-order accurate in time. 
More precisely, if~\eqref{LieProblem} is defined on a finite-dimensional space, and if the operators $A_i$ are smooth enough, 
then $\| \phi(t^n)-\phi^n \| = O(\Delta t)$~\cite{glowinski2003finite}.
In our case, operator $A$ that is associated with problem~\eqref{structure1}-\eqref{dynamic2} 
will be split into a sum of three operators: 
\begin{enumerate}
 \item The time-dependent Stokes problem with suitable boundary conditions involving structure velocity and fluid stress at the boundary.
 \item The fluid advection problem.
 \item The elastodynamics problem for the structure loaded by the fluid pressure.
\end{enumerate}
These sub-problems are coupled via the kinematic coupling condition and via fluid pressure appearing in the elastodynamics problem.
The kinematic coupling condition also plays a key role in writing 
problem~\eqref{structure1}-\eqref{dynamic2} as a first-order system, based on which the Lie splitting can be performed.

\subsection{The first-order system in ALE framework}
\begin{figure}[ht]
 \centering{
 \includegraphics[scale=0.35]{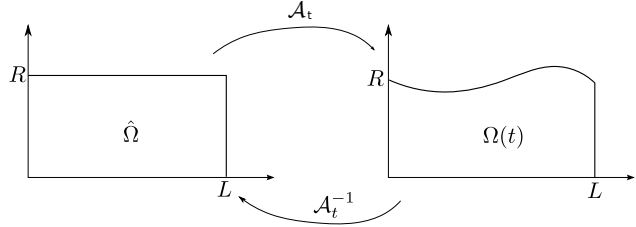}
 }
 \caption{$\mathcal{A}_t$ maps the reference domain $\hat\Omega$ into the current domain $\Omega(t)$. }
\label{ale703}
 \end{figure}

To deal with the motion of the fluid domain we adopt the Arbitrary Lagrangian-Eulerian (ALE) approach~\cite{hughes1981lagrangian,donea1983arbitrary,nobile2001numerical}.
In the context of finite element method approximation of moving-boundary problems,
ALE method deals efficiently 
with the deformation of the mesh, especially near the interface
between the fluid and the structure,
and with the issues 
related to the approximation of the time-derivatives
$\partial \mathbf{u}/\partial t \approx (\mathbf{u}(t^{n+1})-\mathbf{u}(t^{n}))/\Delta t$ which,
due to the fact that $\Omega(t)$ depends on time, is not well defined since
the values $\mathbf{u}(t^{n+1})$ and $\mathbf{u}(t^{n})$ 
correspond to the values of $\mathbf{u}$ defined at two different domains.
ALE approach is based on introducing
a family of (arbitrary, invertible, smooth) mappings ${\cal{A}}_t$ defined on 
a single, fixed, reference domain $\hat\Omega$ such that, for each 
$t \in (t_0, T)$, ${\cal{A}}_t$ maps the reference domain $\hat{\Omega}=(0,L) \times (0,R)$ into the current domain $\Omega(t)$ (see Figure~\ref{ale703}):
$$\mathcal{A}_t : \hat{\Omega} \subset \mathbb{R}^2 \rightarrow \Omega(t) \subset \mathbb{R}^2,  \quad \mathbf{x}=\mathcal{A}_t(\hat{\mathbf{x}}) \in \Omega(t), \quad \textrm{for} \; \hat{\mathbf{x}} \in \hat{\Omega}.$$
%Generally, $\mathcal{A}_t$ should be a homeomorphism, i.e. $\mathcal{A}_t$ is invertible with continuous inverse. 
In our approach, we define $\mathcal{A}_t$ to be a harmonic extension of the mapping
$\tilde g$ that maps the boundary of $\hat\Omega$ to the boundary of $\Omega(t)$ for a given time $t$. 
More precisely, in our case $\hat\Omega :=  (0,L)\times(0,R)$, and so $\mathcal{A}_t$ is a harmonic extension of
$$\tilde{g} : \partial \hat\Omega \to \partial \Omega(t)$$ onto the whole domain $\hat\Omega$, for a given $t:$
\begin{eqnarray*}
 \Delta \mathcal{A}_t &=& 0 \quad \rm{in} \; \hat{\Omega}, \\
\mathcal{A}_t |_{\hat{\Gamma}} &=& \tilde{g}, \\
\mathcal{A}_t |_{\partial \hat{\Omega}\backslash \hat{\Gamma}}&=&0.
\end{eqnarray*}

To rewrite system~\eqref{structure1}-\eqref{dynamic2} in the ALE framework
we notice that for a function $f=f(\mathbf{x},t)$ defined on $\Omega(t) \times (0,T)$ 
the corresponding function $\hat{f} := f \circ \mathcal{A}_t$ defined on $\hat{\Omega} \times (0,T)$
is given by
$$\hat{f}(\hat{\mathbf{x}},t) = f(\mathcal{A}_t(\hat{\mathbf{x}}),t).$$
Differentiation with respect to time, after using the chain rule, gives
\begin{equation}
 \frac{\partial f}{\partial t}\bigg|_{\hat{\boldsymbol x}} =  \frac{\partial {f}}{\partial t}+\mathbf{w} \cdot \nabla f,
\end{equation}
where 
%$\hat{\nabla} = \nabla_{\hat{\mathbf{x}}}$ and 
$\mathbf{w}$ denotes domain velocity given by
\begin{equation}
 \mathbf{w}(\mathbf{x},t) = \frac{\partial \mathcal{A}_t(\hat{\boldsymbol x})}{\partial t}.\bigg. \label{w}
\end{equation}

To write system~\eqref{structure2}-\eqref{dynamic2} in first-order form, we utilize the 
kinematic coupling condition~\eqref{kinematic}. Written in the ALE framework, our problem now reads:
Find $\mathbf{u}= (u_z,u_r)$, $\boldsymbol{\eta}=(\eta_z,\eta_r)$, 
with $\hat{\mathbf{u}}(\hat{\mathbf{x}},t) = \mathbf{u}(\mathcal{A}_t(\hat{\mathbf{x}}),t)$ and
$\hat{\mathbf{u}}|_{\hat{\Gamma}}=\hat{\mathbf{u}}(\hat{z},R,t)$, such that
\begin{eqnarray}
& & \rho_f \bigg( \frac{\partial \mathbf{u}}{\partial t}\bigg|_{\hat{\boldsymbol x}}+ (\mathbf{{u}}-\mathbf{w}) \cdot \nabla \mathbf{u} \bigg)  = \nabla \cdot \boldsymbol\sigma,   \quad \quad \textrm{in}\; \Omega(t) \times (0,T), \label{sys1}\\
& &\nabla \cdot \mathbf{u} = 0  \quad \quad \; \textrm{in}\; \Omega(t) \times (0,T), 
\end{eqnarray}
with the kinematic and dynamic coupling conditions holding on $\Gamma(t)$:
\begin{eqnarray}
& &\frac{\partial \boldsymbol \eta}{\partial t}  = \hat{\mathbf{u}}|_{\hat{\Gamma}}   \quad \quad \; \textrm{on}\; (0,L) \times (0,T), \\
& & \rho_s h \frac{\partial (\hat{u}_z|_{\hat{\Gamma}})}{\partial t}-C_2 \frac{\partial \eta_r}{\partial z}-C_3 \frac{\partial^2 \eta_z}{\partial z^2}-D_2 \frac{\partial (\hat{u}_r|_{\hat{\Gamma}})}{\partial z} -D_3 \frac{\partial ^2 (\hat{u}_z|_{\hat{\Gamma}})}{ \partial z^2}  \nonumber \\
& & \quad = - \ \sqrt{\bigg(1+ \frac{\partial \eta_z}{\partial z} \bigg)^2+\bigg(\frac{\partial \eta_r}{\partial z}\bigg)^2} \widehat{\boldsymbol \sigma \mathbf{n}}|_{\Gamma(t)}  \cdot \mathbf{e_z}   \quad  \; \textrm{on}\; (0,L) \times (0,T), \label{ALEdyn1}\\
& &\rho_s h \frac{\partial (\hat{u}_r|_{\hat{\Gamma}})}{\partial t}+C_0 \eta_r -C_1 \frac{\partial^2 \eta_r}{\partial z^2}  +C_2 \frac{\partial \eta_z}{\partial z} + C_4 \frac{\partial^4 \eta_r}{\partial z^4}+D_0 \hat{u}_r|_{\hat{\Gamma}}\nonumber \\ 
& & -D_1 \frac{\partial^2(\hat{u}_r|_{\hat{\Gamma}})}{ \partial z^2} 
\quad +D_2 \frac{\partial  (\hat{u}_z|_{\hat{\Gamma}})}{ \partial z}+D_4 \frac{\partial ^4 (\hat{u}_r|_{\hat{\Gamma}})}{\partial z^4}\nonumber \\
& & =  - \ \sqrt{\bigg(1+ \frac{\partial \eta_z}{\partial z} \bigg)^2+\bigg(\frac{\partial \eta_r}{\partial z}\bigg)^2} \widehat{\boldsymbol \sigma \mathbf{n}}|_{\Gamma(t)}  \cdot \mathbf{e_r}  
 \quad \; \textrm{on}\; (0,L) \times (0,T),\label{ALEdyn2}
\end{eqnarray}
and the following boundary  conditions on $\Gamma_{\rm in} \cup \Gamma_{\rm out} \cup \Gamma_0$:
\begin{eqnarray}
& & \frac{\partial u_z}{\partial r}(z,0,t) =  u_r(z,0,t) = 0 \quad \textrm{on} \; \Gamma_0, \\
& & \mathbf{u}(0,R,t) = \mathbf{u}(L,R,t) = 0, \quad \boldsymbol\eta|_{z=0,L}=\frac{\partial \eta_r}{\partial z}\bigg|_{z=0,L}=0 , \\
& & \mathbf{\sigma n}|_{in}(0, r,t) = -p_{in}(t) \mathbf{n}|_{in}, \\
& & \mathbf{\sigma n}|_{out}(L,r,t) = -p_{out}(t) \mathbf{n}|_{out} \; \textrm{on} \; (0,R) \times (0,T).
\end{eqnarray}
At time $t = 0$ the following initial conditions are prescribed:
\begin{eqnarray}
& & \mathbf{u}|_{t=0} = \mathbf{0}, \quad \boldsymbol\eta|_{t=0} = \mathbf{0}, \quad \frac{\partial \boldsymbol\eta}{\partial t}\bigg|_{t=0} = \mathbf{0}.  \label{sys2}
\end{eqnarray}
Notice how the kinematic coupling condition is used in~\eqref{ALEdyn1} and~\eqref{ALEdyn2} to rewrite the viscous part of the structure equations in terms of the 
trace of the fluid velocity on $\Gamma(t)$. This will be used in the splitting algorithm described below. 
\vskip 0.1in
\noindent
{\bf Remark 2.} As shown in~\cite{duarte2004arbitrary}, if we discretise~\eqref{w} as
\begin{equation}
\boldsymbol w (\boldsymbol x, \tau) = \frac{\mathcal{A}_{\tau}(\hat{\boldsymbol x})- \hat{\boldsymbol x}}{\tau} + O(\tau),
\end{equation}
we obtain a linear affine transformation for $\mathcal{A}_{\tau}$
\begin{equation}
\mathcal{A}_{\tau}(\hat{\boldsymbol x}) = \hat{\boldsymbol x}+\tau \boldsymbol w (\boldsymbol x, \tau) +O(\tau).
\end{equation}
It can be easily shown that, using this transformation, spatial partial derivatives of a function on a domain $\Omega(\tau)$ 
are equal to the derivatives of the same function on the reference domain $\hat{\Omega}$, plus an error $O(\tau)$
\cite{duarte2004arbitrary}. We avoid dealing with this problem by writing only the time-derivative on the reference domain,
and leaving the spatial derivatives evaluated on the current domain. 
\vskip 0.1in

\subsection{Details of the operator-splitting scheme}
We split the first-order system~\eqref{sys1}-\eqref{sys2} into three sub-problems.
The fluid problem will be split into its viscous part and the pure advection part (incorporating the fluid and ALE advection simultaneously).
The fluid stress $\widehat{\mathbf{\boldsymbol\sigma n}}$ will be split into two parts, Part I and Part II:
\begin{equation}\label{beta}
\widehat{\mathbf{\boldsymbol\sigma n}} = \underbrace{\widehat{\mathbf{\boldsymbol\sigma n}}+\beta\widehat{p\mathbf{n}}}_{(I)}
\underbrace{- \beta \widehat{p\mathbf{n}}}_{(II)},
\end{equation}
where $\beta$ is a number between $0$ and $1$, $0 \le \beta \le 1$, 
with $\beta=0$ corresponding to the splitting introduced in~\cite{guidoboni2009stable}.
As will be shown later, the accuracy of the scheme changes as the value of $\beta$ increases from 0 to 1.
The numerical results presented in this manuscript will correspond to the value of $\beta = 1$, since
our numerical investigation showed that $\beta = 1$ provides highest accuracy for Examples 1 and 2 presented 
in this manuscript.

The viscoelastic structure equations will be split into their viscous part and the elastic part.
These are combined into a splitting algorithm in the following three steps.
\begin{itemize}
 \item \textbf{Step 1.} Step 1 involves solving the time-dependent Stokes problem, incorporating
 the viscous part of the structure and Part I of the fluid stress via a Robin-type boundary condition. The time-dependent Stokes problem is solved
 on a fixed domain $\Omega(t^n)$.
 The problem reads as follows:
 
Find $\mathbf{u}, p$ and $\boldsymbol\eta$, with $\hat{\mathbf{u}}(\hat{\mathbf{x}},t) = \mathbf{u}(\mathcal{A}_t(\hat{\mathbf{x}}),t)$
such that for $t\in (t^n, t^{n+1})$, with $p^n$ and $\boldsymbol \eta^n$ obtained at the previous time step:
\begin{equation}\label{step1}
 \left\{\begin{array}{l@{\ }}
 \rho_f \frac{\partial \mathbf{u}}{\partial t}\big|_{\hat{\boldsymbol x}} =\nabla \cdot \boldsymbol{\sigma},   \quad \nabla \cdot \mathbf{u}=0 \quad \textrm{in} \; \Omega(t^n) \\ \\
 \frac{\partial \boldsymbol{\eta}}{\partial t} (\hat{z},t)  = 0 \quad \textrm{on} \; (0,L) \times(t^n, t^{n+1}), \\ \\ 
 \rho_s h \frac{\partial (\hat{u}_z|_{\hat{\Gamma}})}{\partial t}-D_2 \frac{\partial (\hat{u}_r|_{\hat{\Gamma}})}{\partial z} -D_3 \frac{\partial ^2 (\hat{u}_z|_{\hat{\Gamma}})}{ \partial z^2} +\beta \sqrt{(1+\frac{\partial \eta^n_z}{\partial z})^2+(\frac{\partial \eta_r^n}{\partial z})^2} (\widehat{p^n \mathbf{n}^n})|_{\Gamma(t^n)} \cdot \mathbf{e_z}  \\
 \quad \quad \quad= -\sqrt{(1+\frac{\partial \eta^n_z}{\partial z})^2+(\frac{\partial \eta_r^n}{\partial z})^2} ( \widehat{\boldsymbol \sigma \boldsymbol {n}^n})|_{\Gamma(t^n)} \cdot \mathbf{e_z}    \;  \; \textrm{on}\; (0,L)\times(t^n, t^{n+1}), \\ \\
\rho_s h \frac{\partial (\hat{u}_r|_{\hat{\Gamma}})}{\partial t}+D_0 \hat{u}_r|_{\hat{\Gamma}}-D_1 \frac{\partial^2(\hat{u}_r|_{\hat{\Gamma}})}{ \partial z^2} +D_2 \frac{\partial  (\hat{u}_z|_{\hat{\Gamma}})}{ \partial z}+D_4 \frac{\partial ^4 (\hat{u}_r|_{\hat{\Gamma}})}{\partial z^4}\\
 \quad \quad \quad +\beta \sqrt{(1+\frac{\partial \eta^n_z}{\partial z})^2+(\frac{\partial \eta_r^n}{\partial z})^2} (\widehat{p^n \mathbf{n}^n})|_{\Gamma(t^n)} \cdot \mathbf{e_r}  \\
 \quad \quad \quad = -\sqrt{(1+\frac{\partial \eta^n_z}{\partial z})^2+(\frac{\partial \eta_r^n}{\partial z})^2} ( \widehat{\boldsymbol \sigma \boldsymbol {n}^n})|_{\Gamma(t^n)} \cdot \mathbf{e_r}   
 \quad \textrm{on} \; (0,L)\times(t^n, t^{n+1}), \end{array} \right.  \nonumber
\end{equation}
with the following boundary conditions on $\Gamma_{\rm in}\cup\Gamma_{\rm out}\cup\Gamma_0$: 
\begin{equation*}
  \frac{\partial u_z}{\partial r}(z,0,t) = \quad u_r(z,0,t) = 0 \quad \textrm{on} \; \Gamma_0,
\end{equation*}
\begin{equation*}
  \mathbf{u}(0,R,t) = \mathbf{u}(L,R,t) = 0, 
\end{equation*}
\begin{equation*}
 \mathbf{\boldsymbol\sigma n}|_{in} = -p_{in}(t)\mathbf{n}|_{in}\  {\rm on}\ \Gamma_{\rm in}, \; \; \mathbf{\boldsymbol\sigma n}|_{out} = -p_{out}(t)\mathbf{n}|_{out}  \ {\rm on}\ 
 \Gamma_{\rm out},
\end{equation*}
and initial conditions
\begin{equation*}
\mathbf{u}(t^n)=\mathbf{u}^n, \quad \boldsymbol{\eta}(t^n)=\boldsymbol\eta^n.
\end{equation*}
Then set $\mathbf{u}^{n+1/3}=\mathbf{u}(t^{n+1}), \; \boldsymbol\eta^{n+1/3}=\boldsymbol\eta(t^{n+1}),\; p^{n+1}=p(t^{n+1}).$ \\
Note that here we used only Part I of the fluid stress.
\item \textbf{Step 2:} Solve the fluid and ALE advection sub-problem defined on a fixed domain $\Omega(t^n)$. The problem reads:
Find $\mathbf{u}$ and $\boldsymbol {\eta}$ with $\hat{\mathbf{u}}(\hat{\mathbf{x}},t) = \mathbf{u}(\mathcal{A}_t(\hat{\mathbf{x}}),t)$,
such that for $t\in (t^n, t^{n+1})$
\begin{equation*}
 \left\{\begin{array}{l@{\ }} 
 \frac{\partial \mathbf{u}}{\partial t}\big|_{\hat{\boldsymbol x}} + (\mathbf{u}^{n+1/3}-\mathbf{w}^{n+1/3}) \cdot \nabla \mathbf{u}= 0,   \quad \textrm{in} \; \Omega(t^n)  \\ \\
 \frac{\partial \boldsymbol{\eta}}{\partial t} (\hat{z},t)  = 0 \quad \textrm{on} \; (0,L)\times(t^n, t^{n+1}),\\ \\
\rho_s h_s \frac{\partial (\hat{\mathbf{u}}|_{\hat{\Gamma}})}{\partial t} = 0, \quad \textrm{on} \; (0,L)\times(t^n, t^{n+1}), \end{array} \right. 
\end{equation*}
with boundary conditions: 
$$\mathbf{u}=\mathbf{u}^{n+1/3} \  \; \textrm{on} \; \Gamma_{-}^{n+1/3}, \; \textrm{where}$$ 
$$\Gamma_{-}^{n+1/3} = \{\mathbf{x} \in \mathbb{R}^2 | \mathbf{x} \in \partial \Omega(t^n), (\mathbf{u}^{n+1/3}-\mathbf{w}^{n+1/3})\cdot \mathbf{n} <0 \},$$ 
and initial conditions
$$
\mathbf{u}(t^n)=\mathbf{u}^{n+1/3}, \quad \boldsymbol{\eta}(t^n)=\boldsymbol\eta^{n+1/3}.
$$
Then set $\mathbf{u}^{n+2/3}=\mathbf{u}(t^{n+1}), \; \boldsymbol\eta^{n+2/3}=\boldsymbol\eta(t^{n+1}).$
\item \textbf{Step 3:}  Step 3 involves solving the elastodynamics problem for the location of the deformable boundary
by involving the elastic part of the structure which is loaded by Part II of the normal fluid stress.
Additionally, the fluid and structure communicate via the kinematic lateral boundary condition which gives
the velocity of the structure in terms of the trace of the fluid velocity, taken initially to be the value from the previous step.
The problem reads:
Find $\hat{\mathbf{u}}$ and $\boldsymbol {\eta}$, with $p^{n+1}$ computed in Step 1 and $\boldsymbol \eta^n$ obtained at the previous time step, such that for $t\in(t^n,t^{n+1})$ 
\begin{equation*}
 \left\{\begin{array}{l@{\ }} 
\frac{\partial \mathbf{u}}{\partial t}\big|_{\hat{\boldsymbol x}} = 0,   \quad \textrm{in} \; \Omega(t^n) \\ \\
 \frac{\partial \boldsymbol{\eta}}{\partial t} (z,t)  = \hat{\mathbf{u}}|_{\hat{\Gamma}} \quad \textrm{on} \; (0,L)\times(t^n, t^{n+1}),\\ \\
\rho_s h \frac{\partial( \hat{u}_z|_{\hat{\Gamma}})}{\partial t}-C_2 \frac{\partial \eta_r}{\partial z}-C_3 \frac{\partial^2 \eta_z}{\partial z^2} \\
\quad \quad \quad = \beta \sqrt{(1+\frac{\partial \eta^n_z}{\partial z})^2+(\frac{\partial \eta_r^n}{\partial z})^2} ( \widehat{p^{n+1} \mathbf{n}^n})|_{\Gamma(t^n)} \cdot \mathbf{e_z}  \quad \textrm{on} \; (0,L)\times(t^n, t^{n+1}),\\ 
\rho_s h \frac{\partial (\hat{u}_r|_{\hat{\Gamma}})}{\partial t}+C_0 \eta_r -C_1 \frac{\partial^2 \eta_r}{\partial z^2}  +C_2 \frac{\partial \eta_z}{\partial z} + C_4 \frac{\partial^4 \eta_r}{\partial z^4} \\
\quad \quad \quad =  \beta \sqrt{(1+\frac{\partial \eta^n_z}{\partial z})^2+(\frac{\partial \eta_r^n}{\partial z})^2} ( \widehat{p^{n+1} \mathbf{n}^n})|_{\Gamma(t^n)} \cdot \mathbf{e_r}  \quad \textrm{on} \; (0,L)\times(t^n, t^{n+1}), \end{array} \right.  
\end{equation*}
with boundary conditions: 
\begin{equation*}
  \boldsymbol\eta|_{z=0,L}=\frac{\partial \eta_r}{\partial z}|_{z=0,L}=0; 
\end{equation*}
and initial conditions:
$$
\mathbf{u}(t^n)=\mathbf{u}^{n+2/3}, \quad \boldsymbol{\eta}(t^n)=\boldsymbol\eta^{n+2/3}.
$$
Then set $\mathbf{u}^{n+1}=\mathbf{u}(t^{n+1}), \; \boldsymbol{\eta}^{n+1}=\boldsymbol{\eta}(t^{n+1}).$ \\
Do $t^n=t^{n+1}$ and return to Step 1.
\end{itemize}

\if 1 = 0
Once we calculate the displacement $\boldsymbol\eta$, we can define the ALE mapping as an extension of $\boldsymbol\eta$ onto $\hat\Omega$:
\begin{equation}
 \mathcal{A}_t(\mathbf{x_0})=\mathbf{x_0} + \textrm{Ext}(\boldsymbol\eta(\mathbf{x_0},t)).
\end{equation}
A common choice for the operator Ext, also used here, is a harmonic extension in the reference domain. Now, the domain velocity $\mathbf{\hat{w}}$ can be easily calculated as
\begin{equation}
\mathbf{\hat{w}}(\mathbf{x_0},t) = \textrm{Ext}(\partial_t \boldsymbol\eta(\mathbf{x_0},t)).
\end{equation}
\fi

\vskip 0.1in
\noindent
{\bf Remark 3.} Note that the outward normal to the lateral boundary can be written as
\begin{equation}
\boldsymbol n = \frac{(- \eta_r', 1+\eta_z')}{\sqrt{(\eta_r')^2+(1+\eta_z')^2}}.
\end{equation} 
Using this equality, we can take $\boldsymbol n$ in Step 3 implicitly, which upon substituting $\hat{ \boldsymbol u}|_{\hat{\Gamma}}$ by $\frac{\partial \boldsymbol \eta}{\partial t}$ leads to the following system
\begin{equation*}
 \left\{\begin{array}{l@{\ }} 
\rho_s h \frac{\partial^2 \eta_z}{\partial t^2}-(C_2-\beta \hat{p}|_{\Gamma(t^n)}^{n+1}) \frac{\partial \eta_r}{\partial z}-C_3 \frac{\partial^2 \eta_z}{\partial z^2}=0  \quad \textrm{on} \; (0,L)\times(t^n, t^{n+1}),\\ 
\rho_s h \frac{\partial^2 \eta_r}{\partial t^2}+C_0 \eta_r -C_1 \frac{\partial^2 \eta_r}{\partial z^2}  +(C_2-\beta \hat{p}|_{\Gamma(t^n)}^{n+1}) \frac{\partial \eta_z}{\partial z} + \frac{h^2}{12}C_3 \frac{\partial^4 \eta_r}{\partial z^4} \\
\quad \quad \quad =  \beta \hat{p}|_{\Gamma(t^n)}^{n+1}  \quad \textrm{on} \; (0,L)\times(t^n, t^{n+1}), \end{array} \right.  
\end{equation*}
where $p^{n+1}$ is pressure computed in Step 1.

\vskip 0.1in
\noindent
{\bf Remark 4.} The trace of the pressure, used in Step 3 to load the structure, needs to be well-defined. In general,
one expects the pressure for a Dirichlet problem defined on a Lipschitz domain to be in $L^2(\Omega)$, which is 
not sufficient.
Several works, see e.g., \cite{avalos2008higher,kukavica2012solutions,surulescu2010time}, indicate that, under certain compatibility conditions,
the solution of the related class of moving boundary problems has higher regularity, allowing the definition of the trace of the pressure on the moving boundary.
In fact, we can show that for our problem, under certain compatibility conditions at the corners of the domain, the pressure belongs to $W^{1,8/7}(\Omega)$,
which is more than sufficient for the trace to be well-defined on $\Gamma(t^n)$ \cite{PressureTraceInPreparation}.

\section{Numerical results}\label{sec5}
We present three examples that show the behavior of our scheme for different parameter values.
Example 1 below, corresponds to the benchmark problem suggested 
by Formaggia et al. in~\cite{formaggia2001coupling} to study the behavior of FSI scheme for blood flow.
The structure model in this case is a linearly viscoelastic string model capturing only radial displacement,
with the coefficient that describes bending rigidity given in terms of the shear modulus and Timoshenko
correction factor, which is different from the corresponding Koiter shell model coefficient. 
The structural viscosity constant and the Youngs modulus in this model are both relatively small.
%See, e.g., \cite{Annals_viscoelastic} for the physiological values of these parameters. 

In Example 1b we supplemented this model by an equation describing dynamics of longitudinal displacement, obtained
from the Koiter shell model. The corresponding equation for longitudinal displacement is degenerate, as it does not involve
any spatial derivatives as in the equation for radial displacement, and there are no coupling terms
between the radial and longitudinal displacement. 

In Example 2 we consider the full Koiter shell model capturing both radial and longitudinal displacement,
and the coupling between the two.
The coefficients in the model are given by those associated with the derivation of the Koiter shell,
see~\eqref{coeff}. The values of Youngs modulus, Poisson ratio,
and shell thickness are the same as in Examples 1 and 1b, with the structural viscosity coefficients small,
and related to the one in Example 1. 
This example could be used as a benchmark problem for testing numerical methods for FSI in which the 
structure is modeled as a linearly viscoelastic Koiter shell.

The last example concerns simulation of blood flow through a section of the Common Carotid Artery (CAA).
Numerical simulations were compared with experimental data showing very good agreement.

The following values of fluid and structure parameters are used in Examples 1 and 2, listed below:
{{
\begin{center}
\begin{table}[ht!]
{\small{
\begin{tabular}{|l l l l |}
\hline
\textbf{Parameters} & \textbf{Values} & \textbf{Parameters} & \textbf{Values}  \\
\hline
\hline
\textbf{Radius} $R$ (cm)  & $0.5$  & \textbf{Length} $L$ (cm) & $6$  \\
\hline
\textbf{Fluid density} $\rho_f$ (g/cm$^3$)& $1$ &\textbf{Dyn. viscosity} $\mu$ (poise) & $0.035$    \\
\hline
\textbf{Wall density} $\rho_s $(g/cm$^3$) & $1.1$  & \textbf{Wall thickness} $h_s$ (cm) & $0.1$  \\
\hline
\textbf{Young's mod.} $E $(dynes/cm$^2$) & $0.75 \times 10^6$  & \textbf{Poisson's ratio} $\sigma $ & $0.5$   \\
\hline
\end{tabular}
}}
\caption{Geometric parameters, and fluid and structural parameters that are used in Examples 1 and 2 presented in this section.}
\label{T1}
\end{table}
\end{center}
}}
Parameter $\beta$, introduced in \eqref{beta},  which appears in Step 1 and Step 3 of our numerical scheme, 
can vary between 0 and 1, where the value of $\beta=0$ corresponds to the kinematically coupled scheme presented in~\cite{guidoboni2009stable}. 
The change in $\beta$ is associated with the change in the accuracy of the scheme (not the stability). 
For the test cases presented in Examples 1 and 2, the value of $\beta = 1$ provides the highest accuracy.
We believe that the main reason for the gain 
in accuracy at $\beta = 1$ is the strong coupling between the fluid pressure (which incorporates the leading effect of the fluid loading onto the structure)
and structure elastodynamics, which is established for $\beta = 1$ in Step 3 of the splitting, described above.
Namely, in the case $\beta = 1$, the structure elastodynamics problem is forced by the entire contribution of the fluid pressure,
which appears in Step 3 as a loading onto the sturucture. The elastodynamics of the structural problem is, therefore, directly coupled to the 
entire fluid pressure for $\beta = 1$.
It remains to be investigated how does the optimum choice of $\beta$ (providing highest accuracy for a given problem), depend on the coefficients in the problem.

\subsection{Example 1: The benchmark problem with only radial displacement.}
\begin{figure}[ht!]
 \centering{
 \includegraphics[scale=0.65]{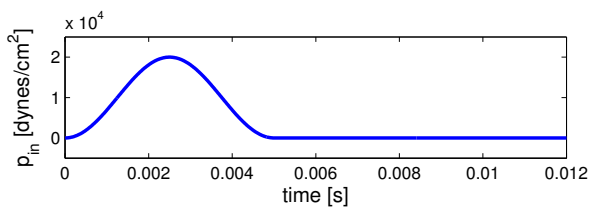}
 }
 \caption{The inlet pressure pulse for Examples 1 and 2. The outlet pressure is kept at 0. }
\label{pressure_pulse}
 \end{figure}

We consider the classical FSI test problem proposed by Formaggia et al.~in~\cite{formaggia2001coupling}.
This problem has been used in several works
as a benchmark problem for testing the results of fluid-structure interaction algorithms for blood flow~\cite{badia2008fluid,nobile2001numerical,badia2008splitting,quaini2009algorithms,guidoboni2009stable}.
The structure model for this benchmark problem is of the form
\begin{equation}\label{nobile}
 \rho_s h \frac{\partial^2 \eta_r}{\partial t^2} - k G h \frac{\partial^2 \eta_r}{\partial z^2} + \frac{E h}{1-\sigma^2} \frac{\eta_r}{R^2} -\gamma \frac{\partial^3 \eta_r}{\partial z^2 \partial t} = f,
\end{equation}
with absorbing boundary conditions at the inlet and outlet boundaries:
\begin{align}
\frac{\partial \eta_r}{\partial t} - \sqrt{\frac{k G}{\rho_s}} \frac{\partial \eta_r}{\partial z} & = 0 \quad \textrm{at} \; z=0 \label{abs1} \\ 
\frac{\partial \eta_r}{\partial t} + \sqrt{\frac{k G}{\rho_s}} \frac{\partial \eta_r}{\partial z} & = 0 \quad \textrm{at} \; z=L. \label{abs2}
\end{align}
Here $G=\frac{E}{2(1+\sigma)}$ is the \textit{shear modulus} and $k$ is the \textit{Timoshenko shear correction factor}.
The flow is driven by the time-dependent pressure data:
\begin{equation}\label{pressure}
 p_{in}(t) = \left\{\begin{array}{l@{\ } l} 
\frac{p_{max}}{2} \big[ 1-\cos\big( \frac{2 \pi t}{t_{max}}\big)\big] & \textrm{if} \; t \le t_{max}\\
0 & \textrm{if} \; t>t_{max}
 \end{array} \right.,   \quad p_{out}(t) = 0 \;\forall t \in (0, T),
\end{equation}
where $p_{max} = 2 \times 10^4$ (dynes/cm$^2$) and $t_{max} = 0.005$ (s). 
The graph of the inlet pressure data versus time, is shown in Figure~\ref{pressure_pulse}.
The values of all the parameters in this model are given in Tables~\ref{T1} and~\ref{T2}.
\begin{center}
\begin{table}[ht]
\begin{center}
\begin{tabular}{| l  l |}
\hline
\textbf{Parameters} & \textbf{Values}  \\
\hline
\hline
\textbf{Shear mod.} $G$(dynes/cm$^2$) & $0.25 \times 10^6$  \\
\hline
 \textbf{Timoshenko  factor} k & 1   \\
\hline
\textbf{Structural viscosity} $\gamma$ (poise cm) & $0.01$   \\
\hline
\end{tabular}
\caption{Example 1: Structural parameters considered in Example 1, in addition to those listed in Table~\ref{T1}.}
\label{T2}
\end{center}
\end{table}
\end{center}
The structural viscosity and Youngs modulus are both very small. For the typical physiological values of
these parameters see, e.g., \cite{Annals_viscoelastic}. 
This means that the arterial wall in this example is rather elastic. The relatively large value of the 
coefficient in front of the second-order derivative with respect to z (describing bending rigidity), 
minimizes the oscillations that 
would normally appear in such structures. See Example 2 for more details.

\begin{figure}[ht!]
 \centering{
 \includegraphics[scale=0.6]{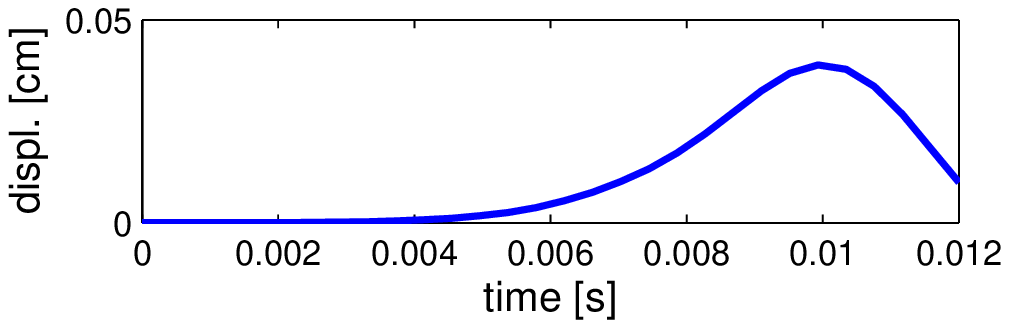}
  \includegraphics[scale=0.6]{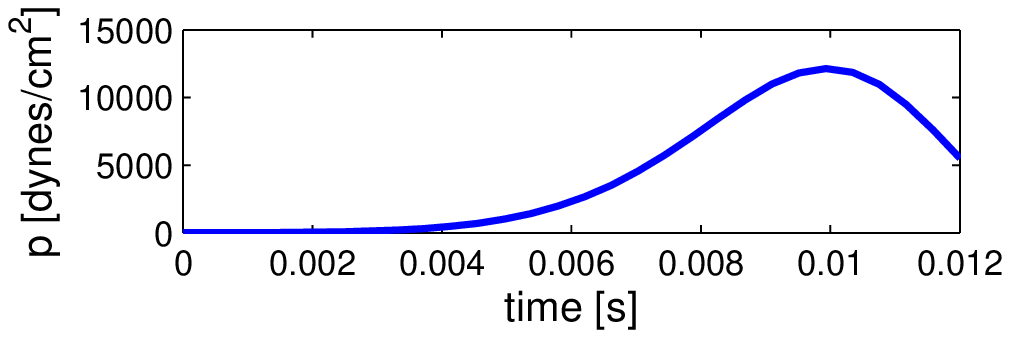}
   \includegraphics[scale=0.6]{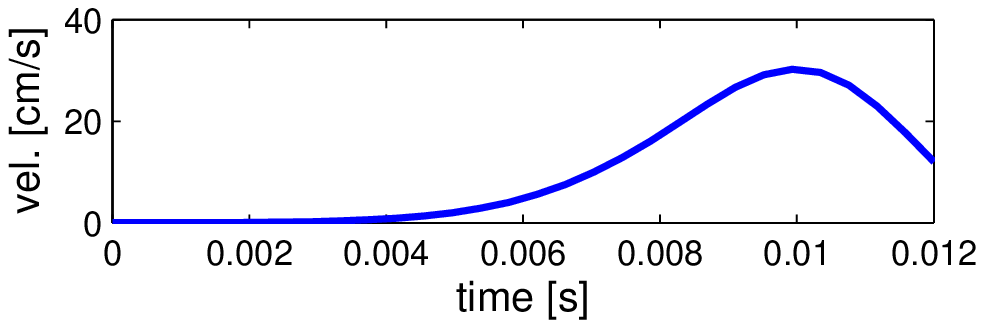}
 }
 \caption{Example 1: Propagation of the inlet pressure pulse in terms of displacement, pressure, and velocity profiles versus time, evaluated at the mid-point of the tube.}
\label{Ex1_propagation}
 \end{figure}

\begin{figure}[ht!]
 \centering{
 \includegraphics[scale=1]{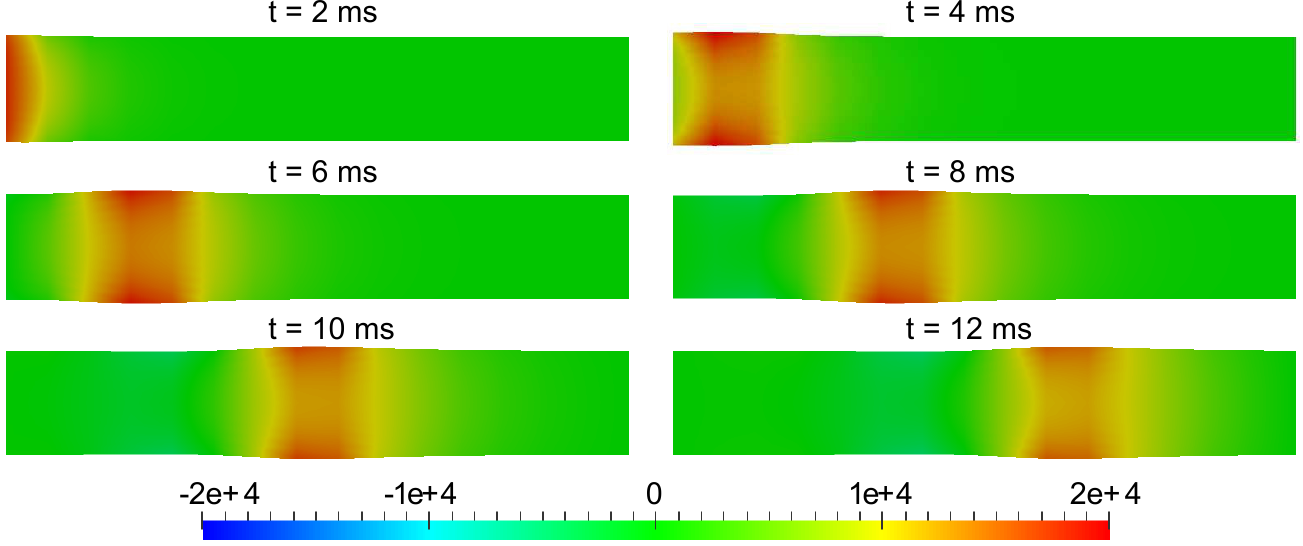}
 }
 \caption{Example 1: Propagation of the pressure wave.}
\label{pressure2D}
 \end{figure}

We implemented this problem in our solver.
Since this model captures only radial displacement the solver was modified accordingly.
The model equation~\eqref{nobile} can be recovered from the Koiter shell model~\eqref{structure2}
by setting the following values for the coefficients:
$$\begin{array}{rlrlrlrlrl}
 C_0 &= \frac{ E h}{R^2(1-\sigma^2)}, \; & C_1 &=k G h, \; & C_2 &=0, \; & C_3 &=0, \; & C_4 &= 0,\\
D_0 &= 0, \; & D_1 &= \gamma, \; & D_2 &= 0, \; & D_3 &= 0, \; & D_4 &= 0.
\end{array}$$
The numerical values of these constants are given in Table~\ref{T3}.
\begin{table}[ht!]
\begin{center}
\begin{tabular}{| l  l  l  l |}
\hline
$C_0 = 4\times 10^5$ & $C_1=2.5 \times 10^4$ & $C_2 = 0 $ & $C_3 = 0$  \\
\hline
$D_0 = 0$  & $D_1 = 10^{-2}$  & $D_2 = 0$ & $D_3 = 0$  \\
\hline
\end{tabular}
\end{center}
\caption{Koter shell model coefficients for Example 1.}
\label{T3}
\end{table}
Homogeneous Dirichlet boundary conditions for the structure in Step 3 were replaced 
with absorbing boundary conditions~\eqref{abs1}-\eqref{abs2}.  
The problem was solved over the time interval $[0,0.012]$s.
\begin{figure}[ht!]
 \centering{
 \includegraphics[scale=0.85]{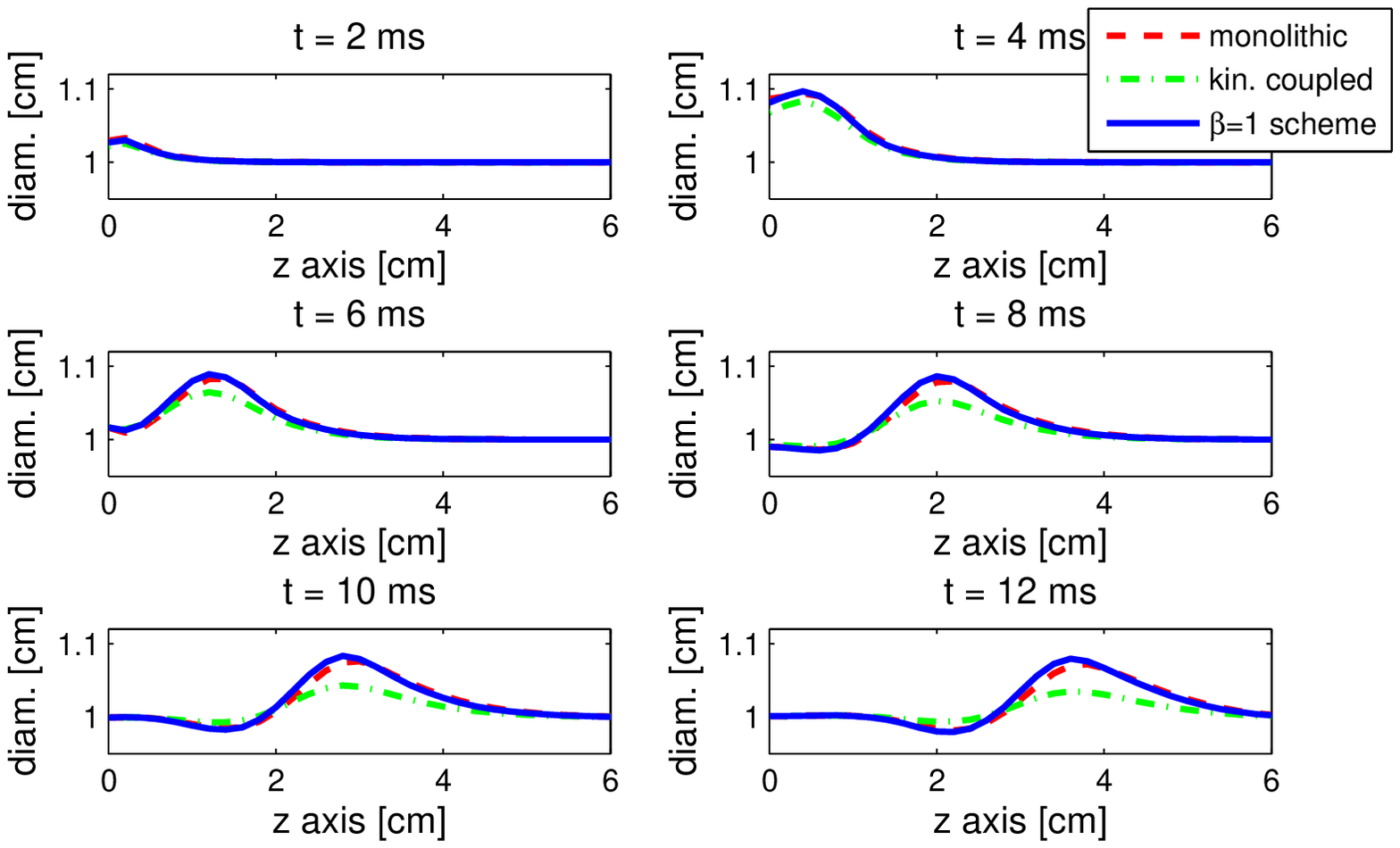}
 }
 \caption{Example 1: Diameter of the tube computed with the kinematically coupled scheme with time step $\triangle t = 10^{-4}$ (dash-dot line), implicit scheme used by Quaini in~\cite{quaini2009algorithms} with the time step $\triangle t = 10^{-4}$ (dashed line) and our scheme with the time step $\triangle t = 10^{-4}$ (solid line). }
\label{displfig}
 \end{figure}
\begin{figure}[ht!]
 \centering{
 \includegraphics[scale=0.85]{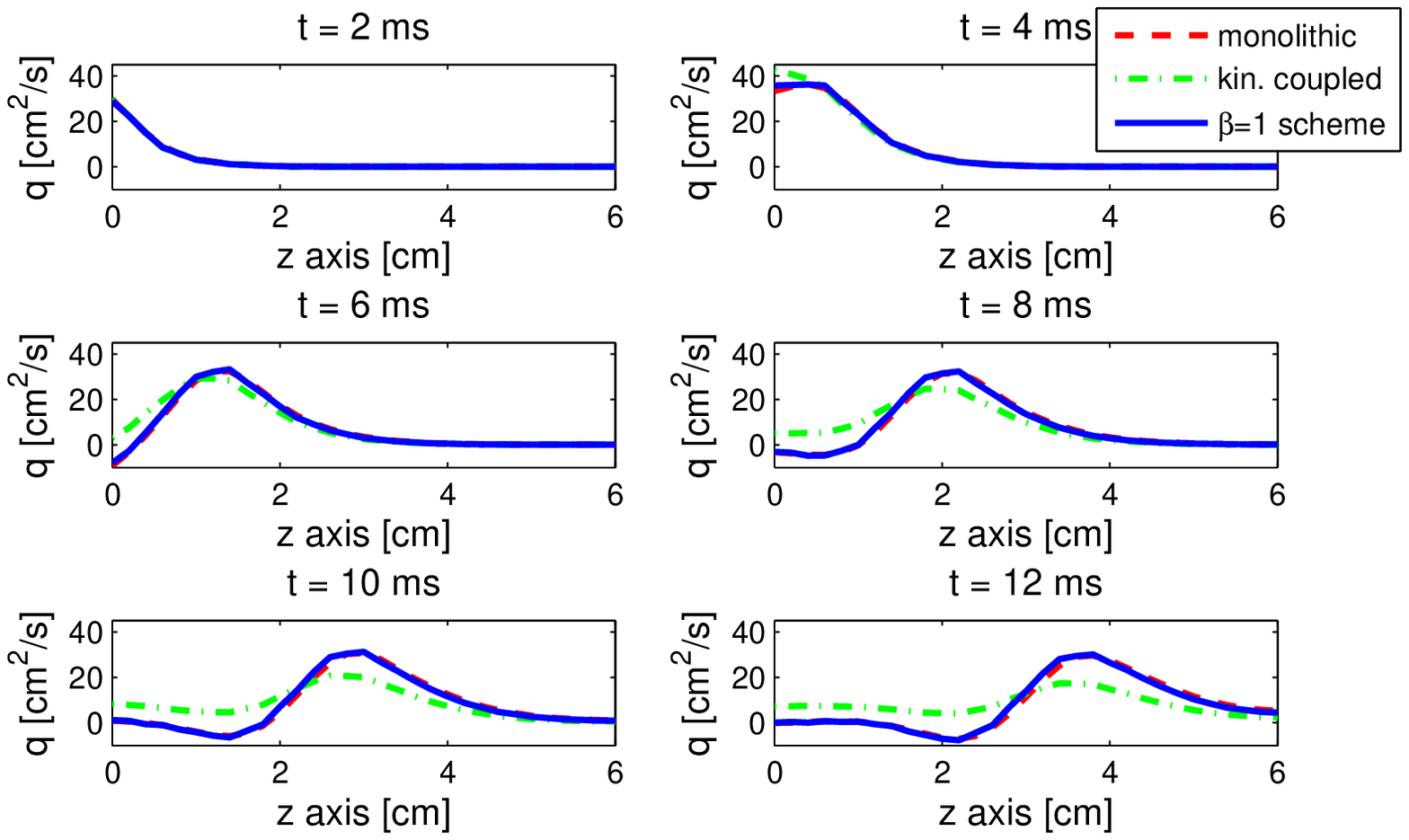}
 }
 \caption{Example 1: Flowrate computed with the kinematically coupled scheme with time step $\triangle t =  10^{-4}$ (dash-dot line), implicit scheme used by Quaini in~\cite{quaini2009algorithms} with the time step $\triangle t = 10^{-4}$ (dashed line) and our scheme with the time step $\triangle t = 10^{-4}$ (solid line). }
\label{flow}
 \end{figure}
\begin{figure}[ht!]
 \centering{
 \includegraphics[scale=0.85]{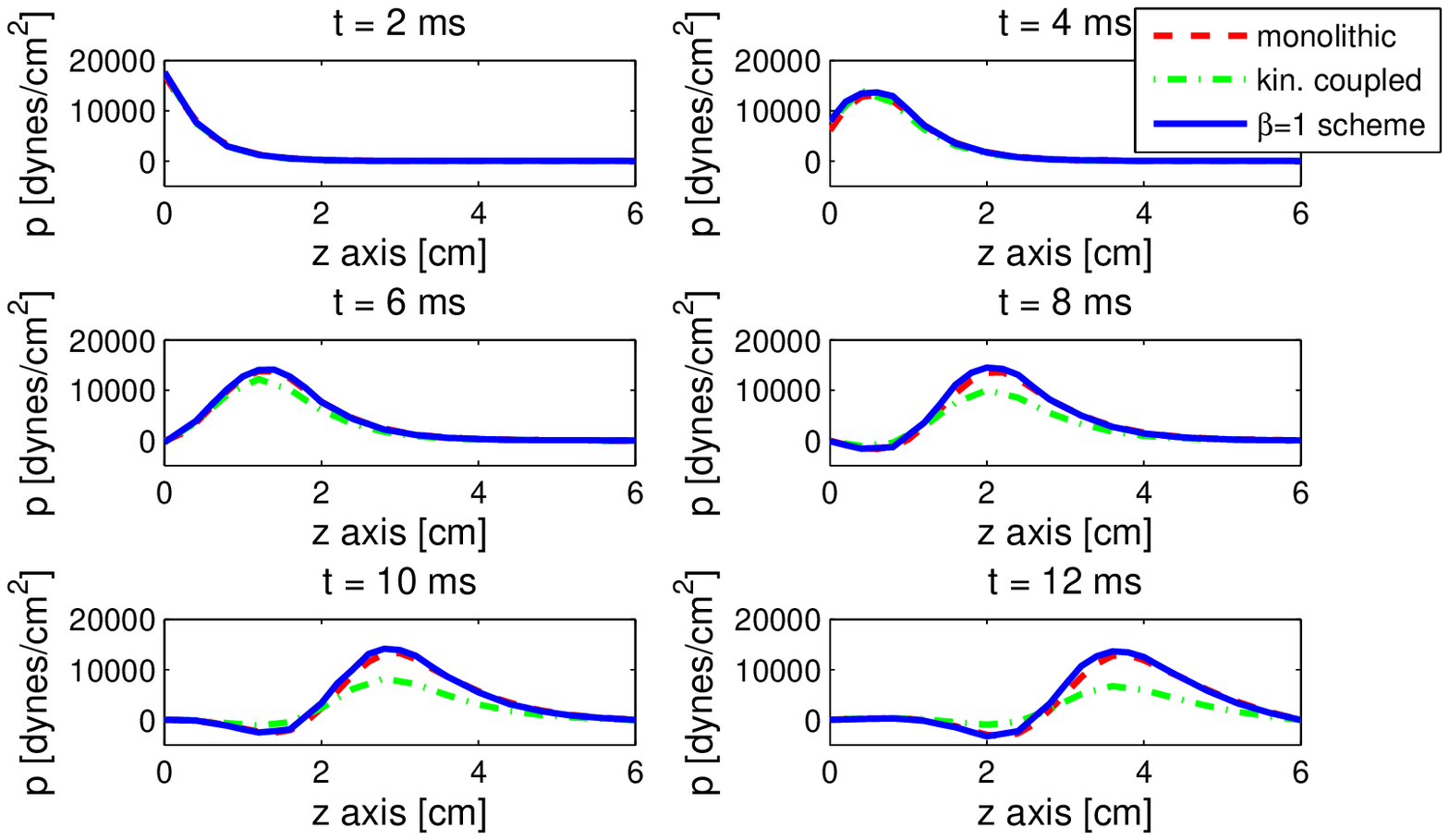}
 }
 \caption{Example 1: Mean pressure computed with the kinematically coupled scheme with time step $\triangle t = 10^{-4}$ (dash-dot line), implicit scheme used by Quaini in~\cite{quaini2009algorithms} with the time step $\triangle t = 10^{-4}$ (dashed line) and our scheme with the time step $\triangle t = 10^{-4}$ (solid line). }
\label{press}
 \end{figure}
Propagation of the inlet pressure pulse in terms of velocity, displacement, and pressure, vs. time, calculated at the mid-point of the tube, are shown in Figure~\ref{Ex1_propagation}.
A 2D cartoon of the pressure pulse propagating in the tube, is shown in Figure~\ref{pressure2D}.

The numerical results obtained using our modification of the kinematically (loosely) coupled scheme proposed in this manuscript,
were compared with the numerical results obtained using the classical kinematically coupled scheme proposed in~\cite{guidoboni2009stable},
and the monolithic scheme proposed in~\cite{quaini2009algorithms}.
Figures~\ref{displfig},~\ref{flow} and~\ref{press} show the comparison between tube diameter, flowrate and mean pressure, respectively, at six different times. 

These results were obtained on the same mesh as the one used for the monolithic scheme in~\cite{quaini2009algorithms}, containing $31 \times 11$ $ \mathbb{P}_1$ fluid nodes. 
More preciesely, we used an isoparametric version of the Bercovier-Pironneau element spaces, also known as 
$ \mathbb{P}_1$-iso-$ \mathbb{P}_2$ approximation in which a coarse mesh is used for the pressure (mesh size $h_p$) and 
a fine mesh for velocity (mesh step $h_v=h_p/2$).

The time step used was $\triangle t = 10^{-4}$ which is the same as the time step used for the monolithic scheme,
and the kinematically coupled scheme.
Due to the splitting error,
it is well-known that classical splitting schemes usually require smaller time step to achieve accuracy comparable to monolithic schemes.
However, the new splitting proposed in this manuscript allows us to use the same time step as in the monolithic method,
obtaining comparable accuracy, as it will be shown next. 
This is exciting since we obtain the same accuracy while
retaining the main benefits of the partitioned schemes, such as 
modularity, simple implementation, and low computational costs.

\begin{figure}[ht!]
 \centering{
 \includegraphics[scale=0.65]{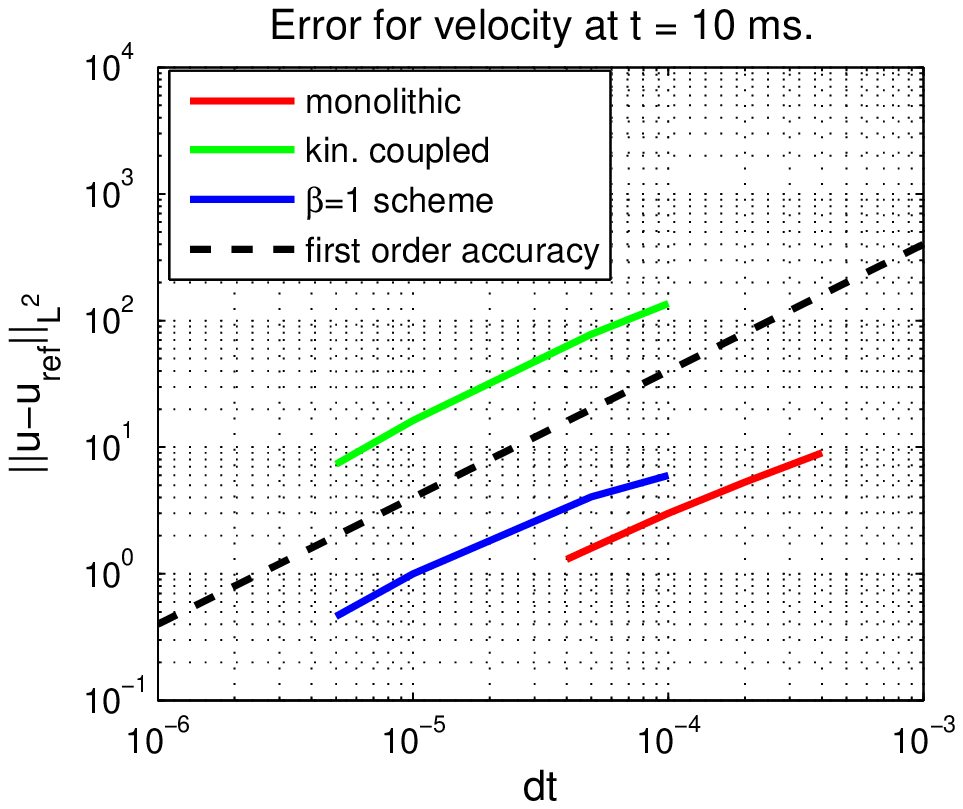}
 \includegraphics[scale=0.65]{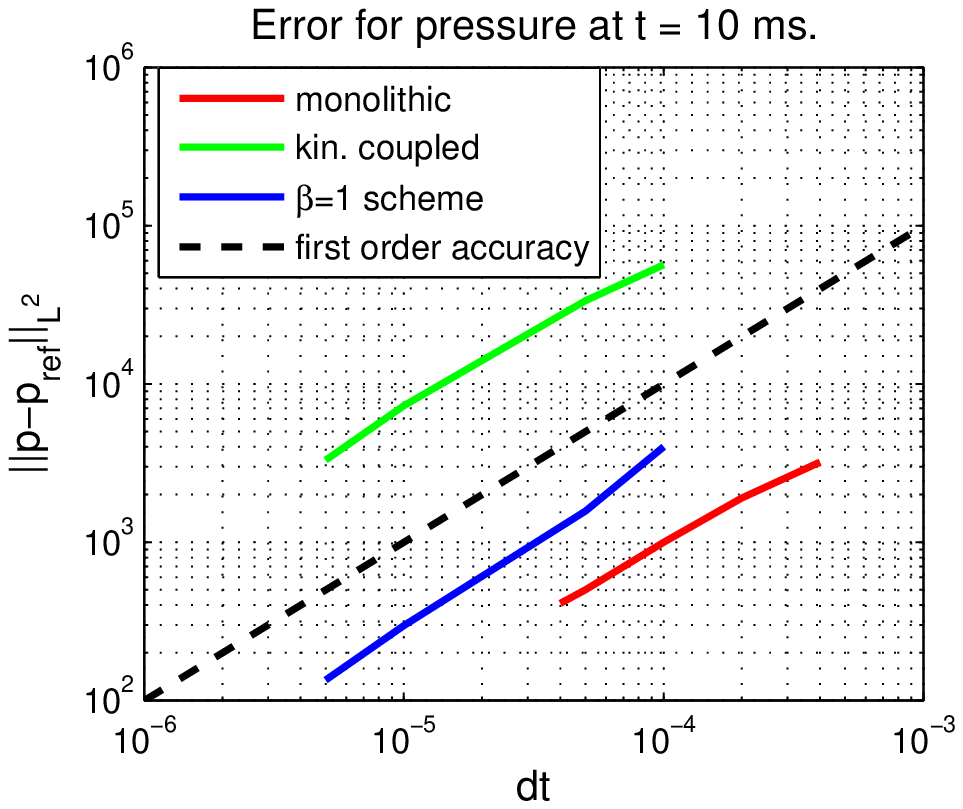} 
 \includegraphics[scale=0.65]{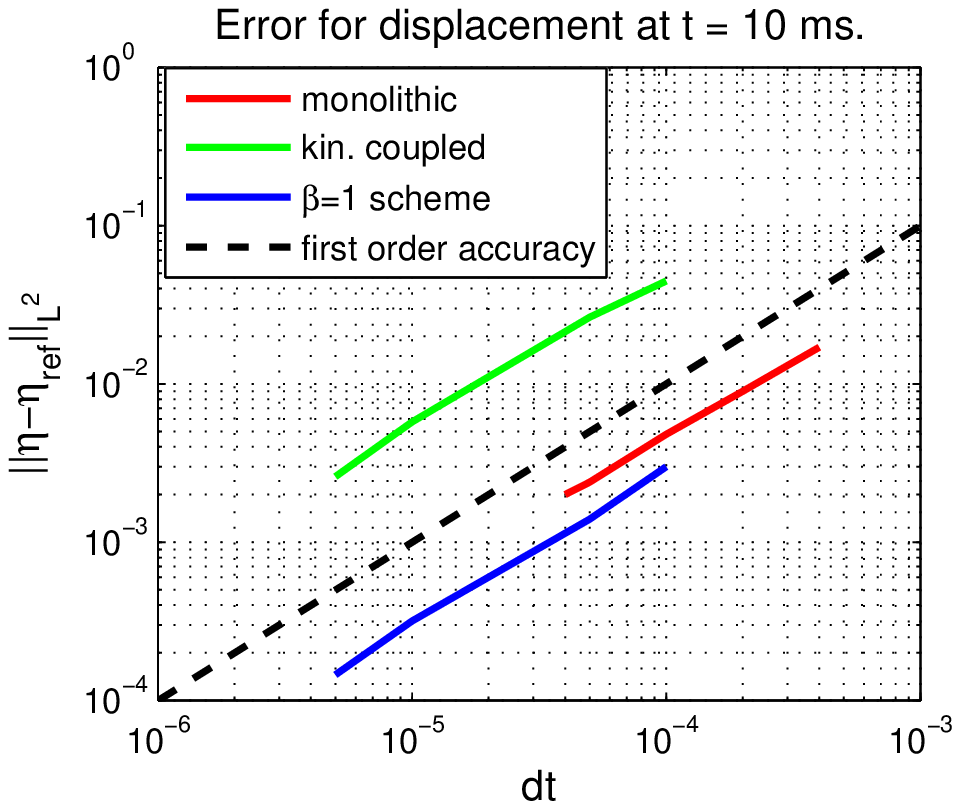}
 }
 \caption{Example 1: Figures show relative errors compared with the kinematically coupled scheme which is first-order accurate in time. Top left: Relative error for fluid velocity at t=10 ms. Top right: Relative error for fluid pressure at t=10 ms. Bottom: Relative error for displacement at t=10 ms.}
\label{error}
 \end{figure}

The kinematically coupled scheme was shown numerically to be first-order accurate  in time and second order accurate in space~\cite{guidoboni2009stable}. 
Second-order accuracy in space is retained in the current kinematically coupled $\beta$-scheme, since the same spatial discretization
of the underlying operators $A_i, i = 1, 2, 3,$ was used in the present manuscript as in \cite{guidoboni2009stable}.
 \begin{table}[ht!]
\begin{center}
{\scriptsize{
\begin{tabular}{| l  c  c  c  c  c  c |}
\hline
$ \triangle t $ & $||p-p_{ref}||_{L^2} $ & $L^2$ order & $||\boldsymbol u-\boldsymbol u_{ref}||_{L^2}$  &$L^2$ order & $ ||\boldsymbol \eta - \boldsymbol \eta_{ref}||_{L^2} $ & $L^2$ order \\
\hline
\hline
$10^{-4}$ & $ 4.01 \textrm{e}+03$  & - & $5.97$  & - & $0.003$  & - \\ 
& $(5.65 \textrm{e}+04)$  & - & $ (136.32)  $& - & $(0.0446)$  & - \\  
\hline  
$5 \times 10^{-5}$ & $ 1.57 \textrm{e}+03$  & 1.35 & $4.05$  & 0.56 & $0.0014$  & 1.1 \\ 
& $(3.36 \textrm{e}+04)$  & (0.75) & $( 77.91)$  & (0.80) & $ (0.0264)$  & (0.75)\\ 
\hline   
$10^{-5}$ & $ 296.36$  & 1.04 & $ 1.0$  & 0.87 & $3.17 \textrm{e}-04$  & 0.92 \\ 
 & $( 7.27 \textrm{e}+03)$  &(0.95) & $(16.27)$  &(0.97) & $ (0.00576)$  & $(0.95)$\\ 
\hline   
$5 \times 10^{-6}$ & $ 134.33$  & 1.14 & $0.46$  & 1.12 & $1.45 \textrm{e}-04$  & 1.13 \\ 
 & $(3.3 \textrm{e}+03)$  & (1.14) & $ (7.36)$  & $(1.14)$ & $(0.0026)$  & $(1.14)$ \\ 
\hline 
\end{tabular}
}}
\end{center}
\caption{Example 1: Convergence in time calculated at $t = 10$ ms.
The numbers in the parenthesis show the convergence rate for the kinematically coupled scheme presented in~\cite{guidoboni2009stable}. }
\label{T4}
\end{table}
However, due to the new time-splitting, the accuracy in time has changed.
Indeed, here we show that this is the case by studying time-convergence of our scheme. 
Figure~\ref{error} shows a comparison between the time convergence of our scheme, the kinematically coupled scheme, and the monolithic scheme used in~\cite{quaini2009algorithms}. 
The reference solution was defined to be the one obtained with $\triangle t = 10^{-6}$. 
We calculated the absolute $L^2$ error for the velocity, pressure and displacement between the reference solution and the solutions obtained using 
$\triangle t = 5\times 10^{-6}, 10^{-5},\; 5 \times 10^{-5}$ and $10^{-4}$.
Figure~\ref{error} shows first-order in time convergence for the velocity, pressure, and displacement obtained by the kinematically coupled scheme, monolithic scheme,
and our scheme. The error of our method is noticably smaller
than the error obtained using the classical kinematically coupled scheme, and is comparable to the error obtained by the monolithic scheme. 
The values of the convergence rates for pressure, velocity, and displacement, calculated using the kinematically coupled schemes, are shown in Table~\ref{T4}.

\subsubsection{Homogeneous Dirichlet vs. absorbing boundary conditions}
\begin{figure}[ht]
 \centering{
 \includegraphics[scale=1.05]{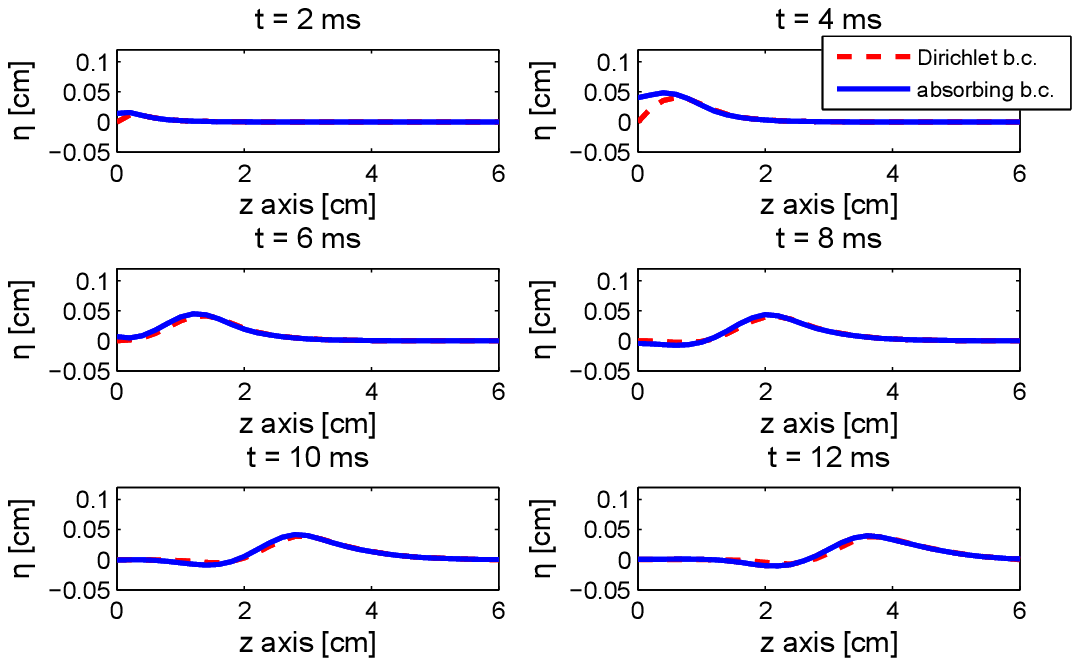}
 }
 \caption{Example 1: Displacement of the structure in the case of absorbing boundary conditions (solid line) and homogeneous Dirichlet boundary conditions (dashed line).}
\label{dir-vs-neu}
 \end{figure}

We give a short remark related to the impact of the homogeneous Dirichlet vs.~absorbing boundary conditions.
Although absorbing boundary conditions for the structure are more realistic in the blood flow application, 
they will only impact the solution near the boundary, except when reflected waves form in which case 
the influence of the boundary conditions is felt everywhere in the domain.  
It was rigorously proved in~\cite{canic2003effective} that in the case of homogeneous Dirichlet inlet/outlet structure data $\boldsymbol\eta = 0$, 
 a boundary layer forms near the inlet or outlet boundaries of the structure to accommodate the transition from 
zero displacement to the displacement dictated by the inlet/outlet normal stress flow data. 
It was proved in~\cite{canic2003effective} that this boundary layer decays exponentially fast away from the inlet/outlet boundaries.
Figure~\ref{dir-vs-neu} depicts a  comparison between the displacement obtained using absorbing boundary conditions, and 
the displacement obtained using homogeneous Dirichlet boundary conditions, showing a boundary layer near the 
inlet boundary where the two solutions differ the most. 

\begin{figure}[ht]
 \centering{
 \includegraphics[scale=1.05]{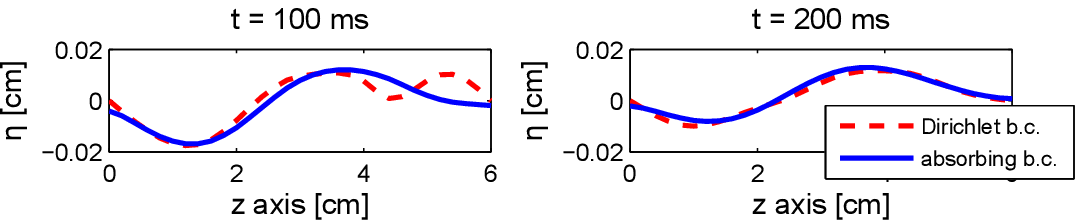}
 }
 \caption{Example 1: Displacement of the structure in the case of absorbing boundary conditions (solid line) and homogeneous Dirichlet boundary conditions (dashed line).}
\label{dirneu200}
 \end{figure}

It is worth pointing out, however, that absorbing boundary conditions help in reducing reflected waves 
that will appear when the propagating wave reaches the outlet boundary and reflects back. 
The ``optimum'' absorbing boundary conditions would have to be designed on the basis of Riemann Invariants (or characteristic variables)
for the hyperbolic problem modeling wave propagation in the structure.
Those conditions, however, are not always easy to calculate, and so approximate 
Riemann Invariant-based absorbing conditions such as \eqref{abs1}-\eqref{abs2} are used.
Figure~\ref{dirneu200} shows the displacement in the case of absorbing boundary conditions \eqref{abs1}, \eqref{abs2},
and homogeneous Dirichlet boundary conditions $\boldsymbol\eta = 0$ at $t=100$ ms and $t=200$ ms. 
Notice how the two solutions differ everywhere in the domain, and how the solution with absorbing boundary conditions 
reduces the amplitude and formation of reflected waves.
%\clearpage

\subsubsection{Example 1b.}
This example is an extension of the benchmark problem by Formaggia et al.~\cite{formaggia2001coupling} studied in Example 1.
The extension concerns inclusion of the longitudinal displacement in the model described in Example 1.

\begin{figure}[ht!]
 \centering{
 \includegraphics[scale=1.05]{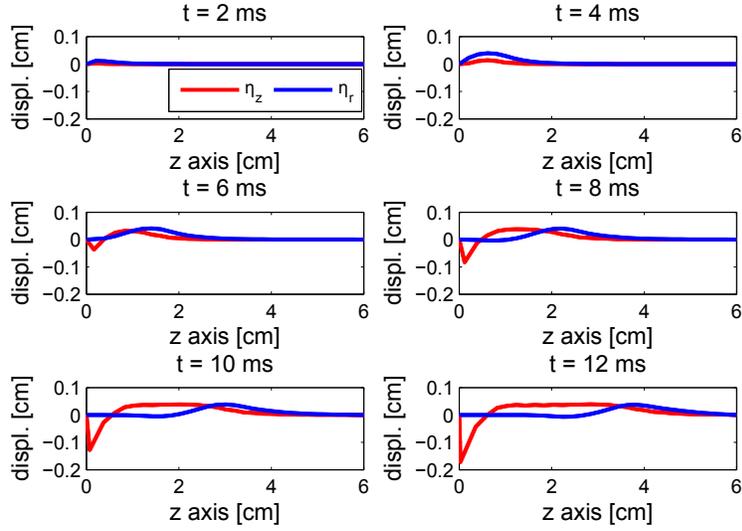}
 }
 \caption{Example 1b: Longitudinal displacement (red) and radial displacement (blue) for the Koiter shell model in Example 2
 calculated with $\Delta t = 10^{-4}$. }
\label{our_displacement}
 \end{figure}
\begin{figure}[ht!]
 \centering{
 \includegraphics[scale=1.05]{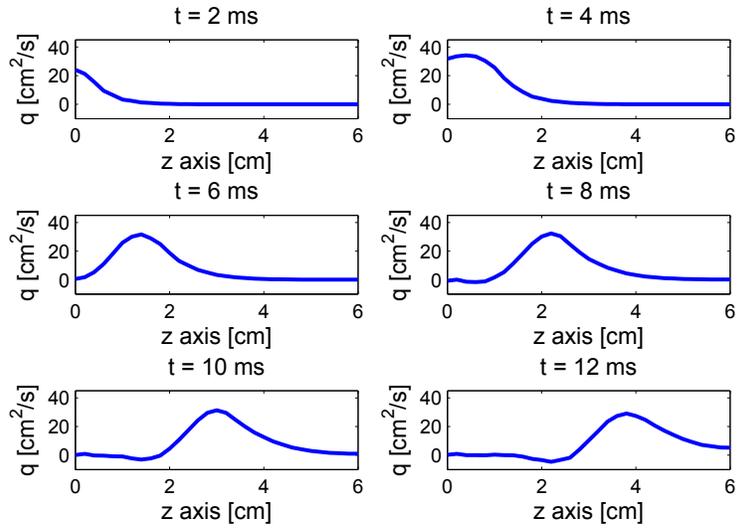}
 }
 \caption{Example 1b: Flow rate for the Koiter shell model in Example 1b
 calculated with $\Delta t = 10^{-4}$.}
\label{our_flow}
 \end{figure}
\begin{figure}[ht!]
 \centering{
 \includegraphics[scale=1.05]{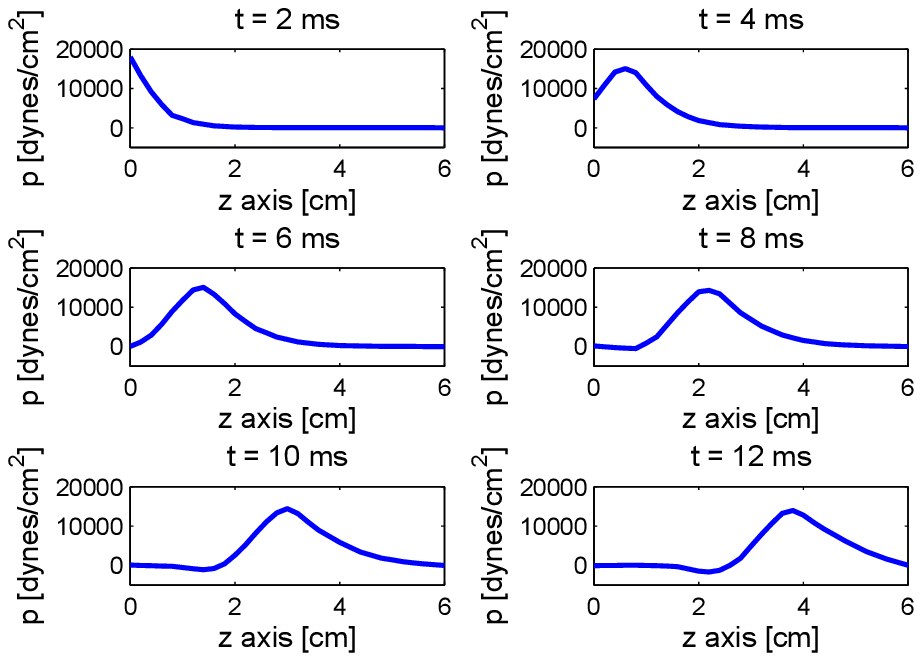}
 }
 \caption{Example 1b:  Mean pressure for the Koiter shell model in Example 1b
 calculated with $\Delta t = 10^{-4}$. }
\label{our_press}
 \end{figure}

Namely, here we explore how the  Koiter shell model \eqref{structure1}-\eqref{structure2}
looks for the coefficients given by those in Example 1.
Namely, the radial displacement satisfies the same model equation as in Formaggia et al.~\cite{formaggia2001coupling},
while the longitudinal displacement satisfies equation \eqref{structure1} in the Koiter shell
model, with  the corresponding coefficients in accordance with equation~\eqref{nobile}.
More precisely, by comparing equations \eqref{structure2} and \eqref{nobile} we observe that only the coefficients
$C_0, C_1$ and $D_1$ are different from zero.
This implies the following Koiter shell model:
\begin{eqnarray}
 \rho_s h \frac{\partial^2 \eta_z}{\partial t^2} &=& f_z,\label{benchmark2}\\
 \rho_s h \frac{\partial^2 \eta_r}{\partial t^2} - k G h \frac{\partial^2 \eta_r}{\partial z^2} + \frac{E h}{1-\sigma^2} \frac{\eta_r}{R^2} 
 -\gamma \frac{\partial^3 \eta_r}{\partial z^2 \partial t} &=& f_r. \label{benchmark1}
\end{eqnarray}
Notice that this problem is degenerate in that the operator associated with the static equilibrium problem is no longer strictly elliptic.
This is due to the fact that the coefficient $C_3$ at the second-order derivative with respect to $z$ of $\eta_z$ in equation \eqref{benchmark2}
is equal to zero.
Nevertheless, we solve the related FSI problem with zero Dirichlet boundary data $\boldsymbol\eta = 0$, and study the
flow driven by the time-dependent pressure data \eqref{pressure} given in Example~1.
The values of the coefficients in the Koiter shell model~\eqref{benchmark1}-\eqref{benchmark2} are equal to those in Example~1.
Figures~\ref{our_displacement},~\ref{our_flow}, and~\ref{our_press} show the displacement, flow rate, and pressure, respectively. 
It is interesting to notice, as is shown in Figure~\ref{our_displacement}, that the magnitude of longitudinal displacement is the same
as the magnitude of radial displacement. 

\subsection{Example 2.} In this test case we consider the full Koiter shell model~\eqref{structure1}-\eqref{structure2}.
This means that all the coefficients $C_i$ and $D_i$ are given in terms of the Young's modulus $E$ and Poisson ratio $\sigma$, and their 
viscous counterparts $E_v$ and $\sigma_v$, through the relationships \eqref{coeff}. Since the terms  containing the 4th and 5th order derivatives are negligible
(it can be shown, using non-dimensional analysis, that these terms are much smaller than the remaining terms), we ignore these terms
in the numerical simulation. We present this example as a benchmark problem for FSI studies in which the structure is modeled by
the Koiter shell model \eqref{structure1}-\eqref{coeff}, which captures both radial and longitudinal displacement.

 The fluid and structure parameters are given
in Tables~\ref{T1} and~\ref{T5}. 
The model and the coefficients are summarized as follows:
\begin{eqnarray*}
\rho_s h \frac{\partial^2 \eta_z}{\partial t^2}-C_2 \frac{\partial \eta_r}{\partial z}-C_3 \frac{\partial^2 \eta_z}{\partial z^2}-D_2 \frac{\partial ^2 \eta_r}{\partial t \partial z} -D_3 \frac{\partial ^3 \eta_z}{\partial t \partial z^2}  = f_z  \label{structure1_ex3}  \\
\rho_s h \frac{\partial^2 \eta_r}{\partial t^2}+C_0 \eta_r -C_1 \frac{\partial^2 \eta_r}{\partial z^2}  +C_2 \frac{\partial \eta_z}{\partial z} +D_0 \frac{\partial \eta_r}{\partial t}-D_1 \frac{\partial^3 \eta_r}{\partial t \partial z^2} 
+D_2 \frac{\partial ^2 \eta_z}{\partial t \partial z}= f_r, \label{structure2_ex3}
\end{eqnarray*}
with
\begin{equation*}
 \begin{array}{rlrlrl}
 C_0 &= \frac{h E}{R^2(1-\sigma^2)}(1+\frac{h^2}{12 R^2}), \; & C_1 &= \frac{h^3}{6} \frac{E \sigma}{R^2 (1-\sigma^2)}, \; & C_2 &=\frac{h}{R}\frac{E \sigma}{1-\sigma^2}, \;   \\
D_0 &= \frac{h}{R^2} C_v(1+\frac{h^2}{12 R^2}), \; & D_1 &= \frac{h^3}{6} \frac{D_v}{R^2}, \; & D_2 &= \frac{h D_v}{R}, \; 
\end{array}
\end{equation*}
where here we take the value of  $D_1$ to be equal to
 $\gamma$ from Examples 1 and 1b, which, from the definition of coefficient $D_1$ above, implies $D_v = 6R^2 \gamma/h^3$. From here we 
determine $C_v = D_v/\sigma_v= 2D_v$.
\begin{table}[ht!]
\begin{center}
\begin{tabular}{|l l|}
\hline
\textbf{Parameters} & \textbf{Values}  \\
\hline
\hline
\textbf{Structural viscosity}   $C_v$(poise cm) & $30$  \\
\hline
 \textbf{Structural viscosity}  $D_v$(poise cm) & $15$  \\
\hline
\end{tabular}
\end{center}
\caption{Structural viscosity parameters for Example 2.}
\label{T5}
\end{table}
The corresponding values of the coefficients are given in Table~\ref{T6}.
\begin{table}[ht!]
\begin{center}
\begin{tabular}{| l  l  l  l|}
\hline
$C_0 = 4.0133 \times 10^5$ & $C_1=333.3$ & $C_2 = 10^5 $  & $C_3 = 10^5$ \\
\hline
$D_0 = 12$  & $D_1 = 10^{-2}$  & $D_2 = 3$ & $D_3 = 3$ \\
\hline
\end{tabular}
\end{center}
\caption{Koter shell model coefficients for Example 2.}
\label{T6}
\end{table}
\begin{figure}[ht!]
 \centering{
\includegraphics[scale=0.95]{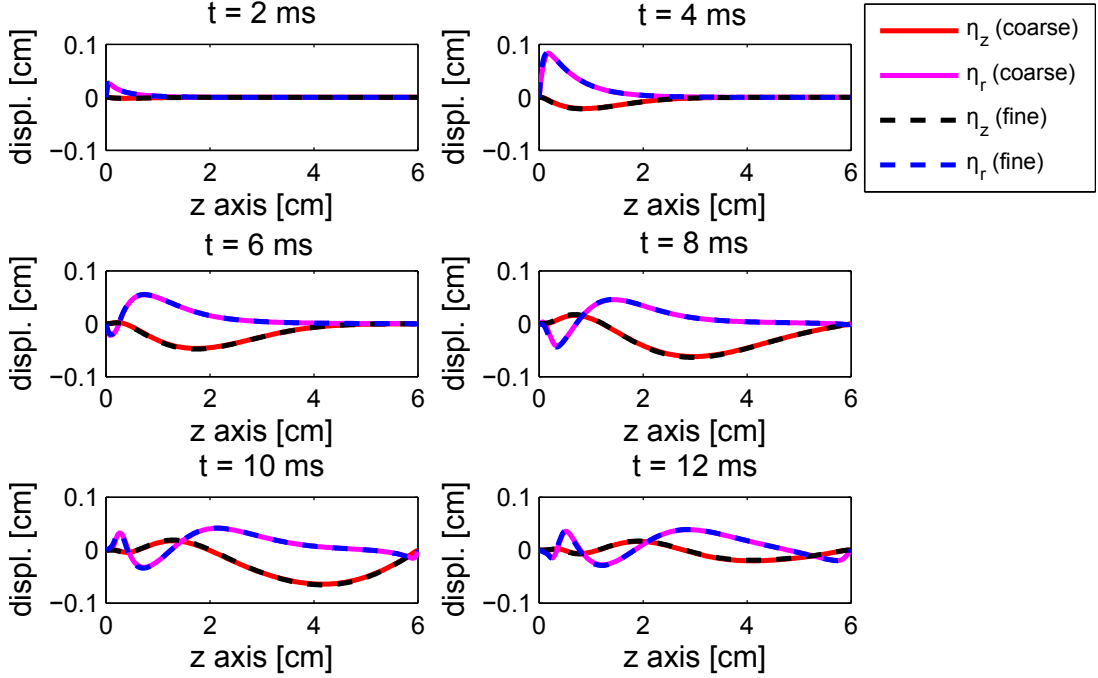}
 }
 \caption{Example 2: Longitudinal displacement $\eta_z$, and radial displacement $\eta_r$ calculated
on a coarse mesh (solid line) and on a fine mesh (dashed line), obtained with $\triangle t = 10^{-4}$.}
\label{displ_elastic}
 \end{figure}
\begin{figure}[ht!]
 \centering{
 \includegraphics[scale=0.95]{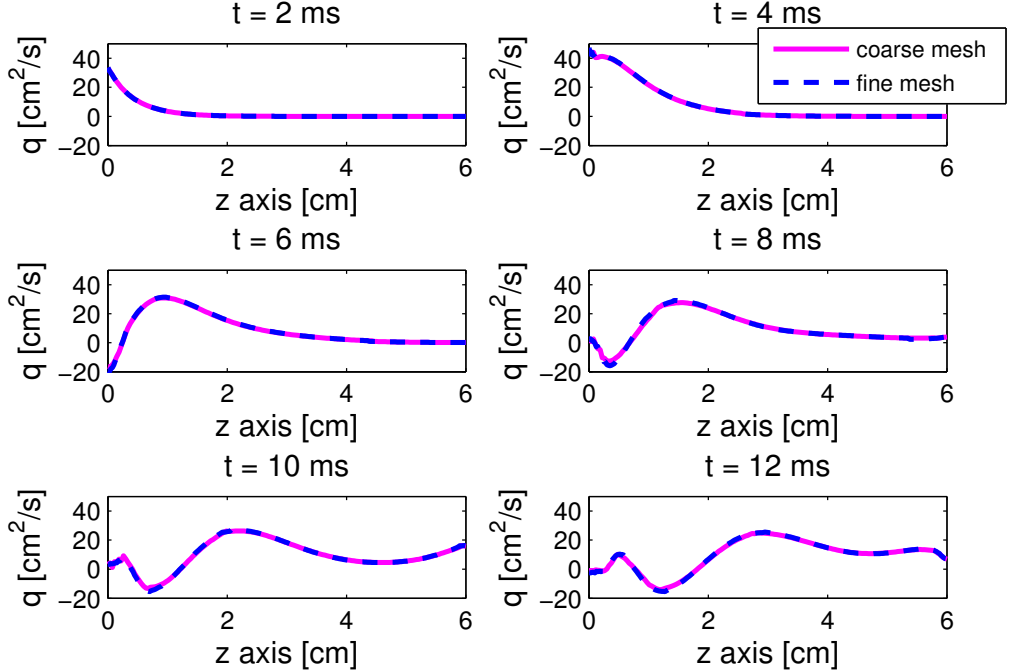}
 }
 \caption{Example 2: Flow rate computed on a coarse mesh (solid line), and on a fine mesh (dashed line),  obtained with $\triangle t = 10^{-4}$.}
\label{flowrate_elastic}
 \end{figure}
\begin{figure}[ht!]
 \centering{
 \includegraphics[scale=0.95]{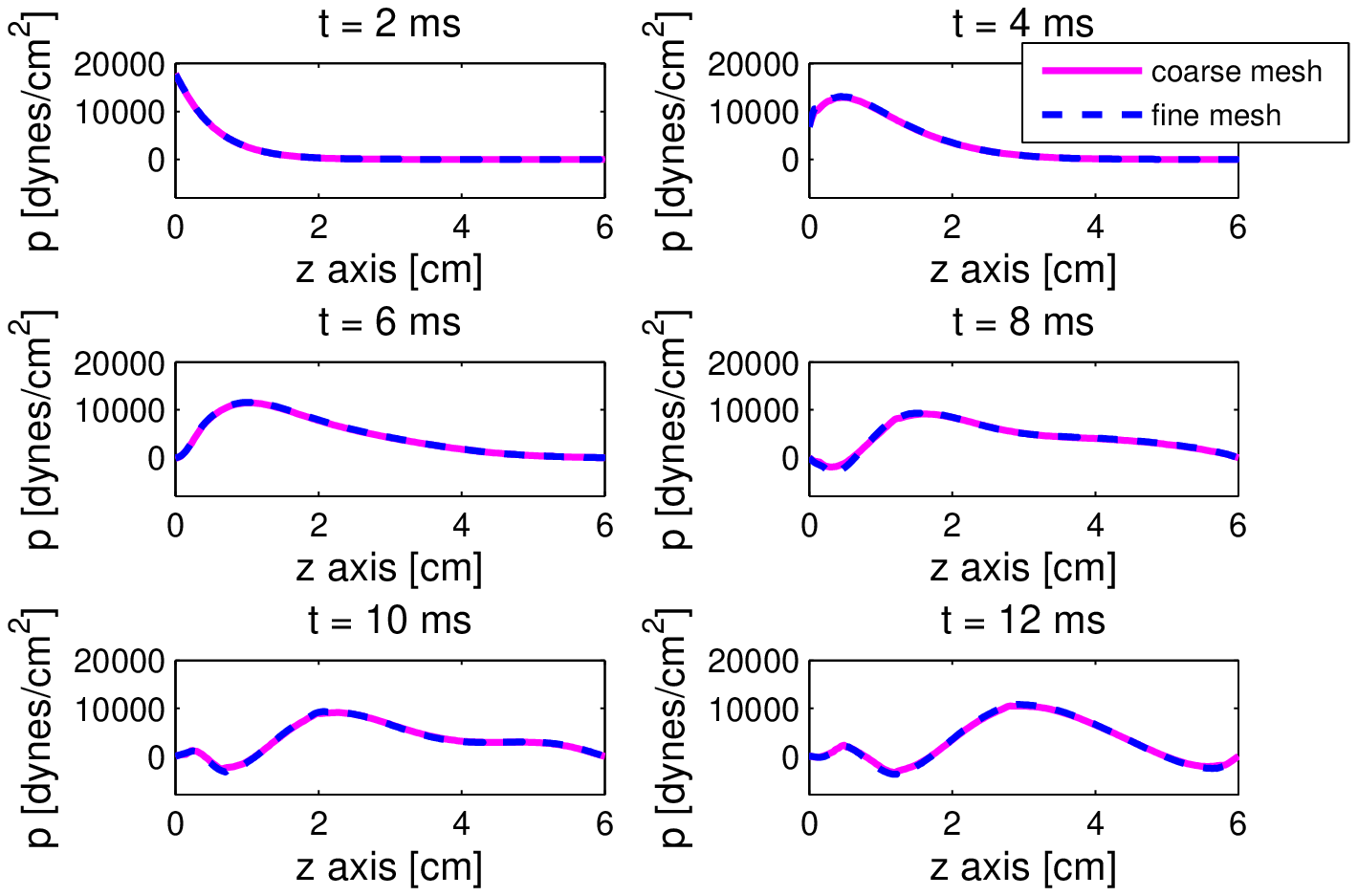}
 }
 \caption{Example 2: Mean pressure computed on a coarse mesh (solid line), and on a fine mesh (dashed line), obtained with $\triangle t = 10^{-4}$.}
\label{pressure_elastic}
 \end{figure}

This model includes the coupling terms between the longitudinal and radial components of the displacement through $C_2\ne0$ and $D_2\ne 0$, 
and the leading-order viscoelastic effects in the radial displacement described by $D_0\ne 0$.
Notice a much smaller value for the coefficient $C_1$ than in Example~1.
Also notice the large coefficient $C_2$ that describes the coupling between the radial and longitudinal 
components of the displacement in the Koiter shell model.  
%Although the viscoelasticity parameters, or the inlet pressure, are not physiologically reasonable, this example can serve as a good benchmark problem to study
%the behavior of FSI schemes for blood flow in which both radial and longitudinal displacement of a thin structure are included.

Figure~\ref{displ_elastic} shows longitudinal and radial displacement of the structure computed on 
two different meshes, evalued at times $t=2,4,6,8,10$ and $12$ ms. 
The coarser mesh is twice as fine as the mesh used in Example 1, so that the triangularization of the coarser mesh
is $h^{coarse}_p = h_p/2$, $h^{coarse}_v = h_v/2$. The fine mesh in this example is twice as fine as the coarse mesh,
namely $h^{fine}_p = h^{coarse}_p/2 = h_p/4$, $h^{fine}_v = h^{coarse}_v/2 = h_v/4$.
Note that the longitudinal displacement is of the same order of magnitude as the radial displacement. 
Figures~\ref{flowrate_elastic} and~\ref{pressure_elastic} show the corresponding flow rate and mean pressure.
The appearance of oscillations near the inlet is physical, and is associated with a very small value of the coefficient $C_1$
when compared to Example 1. An increase in the value of visoelasticity parameters, examined in the next example (physiologically reasonable), 
dampens the oscillations observed in this example.

\begin{figure}[ht!]
 \centering{
 \includegraphics[scale=0.65]{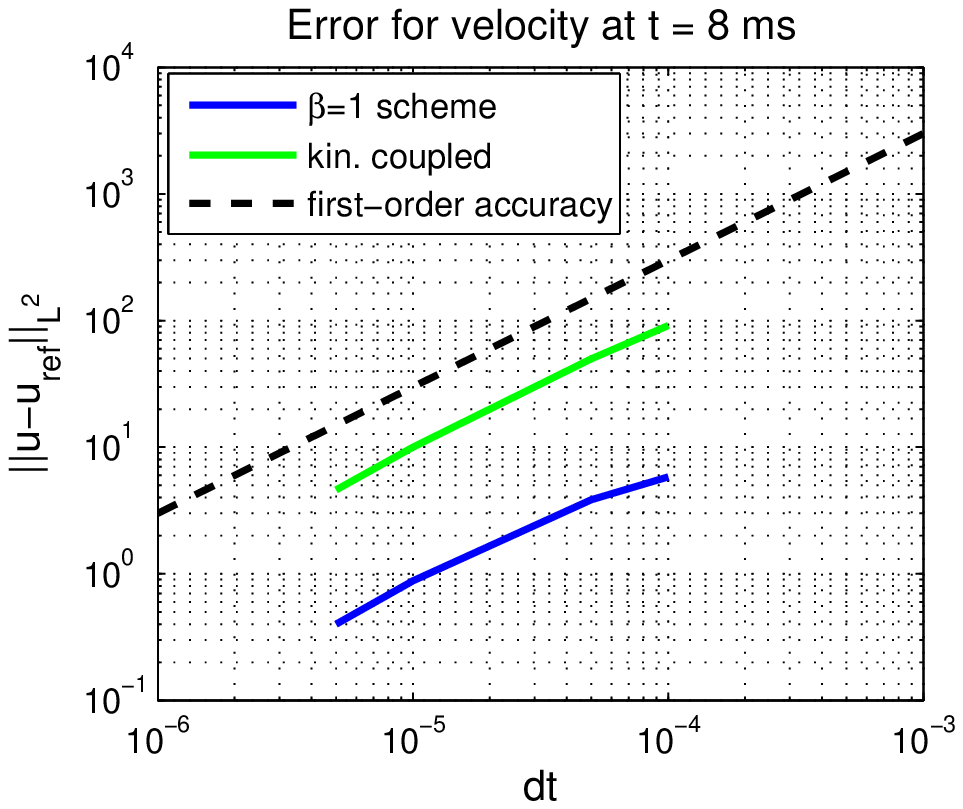}
 \includegraphics[scale=0.65]{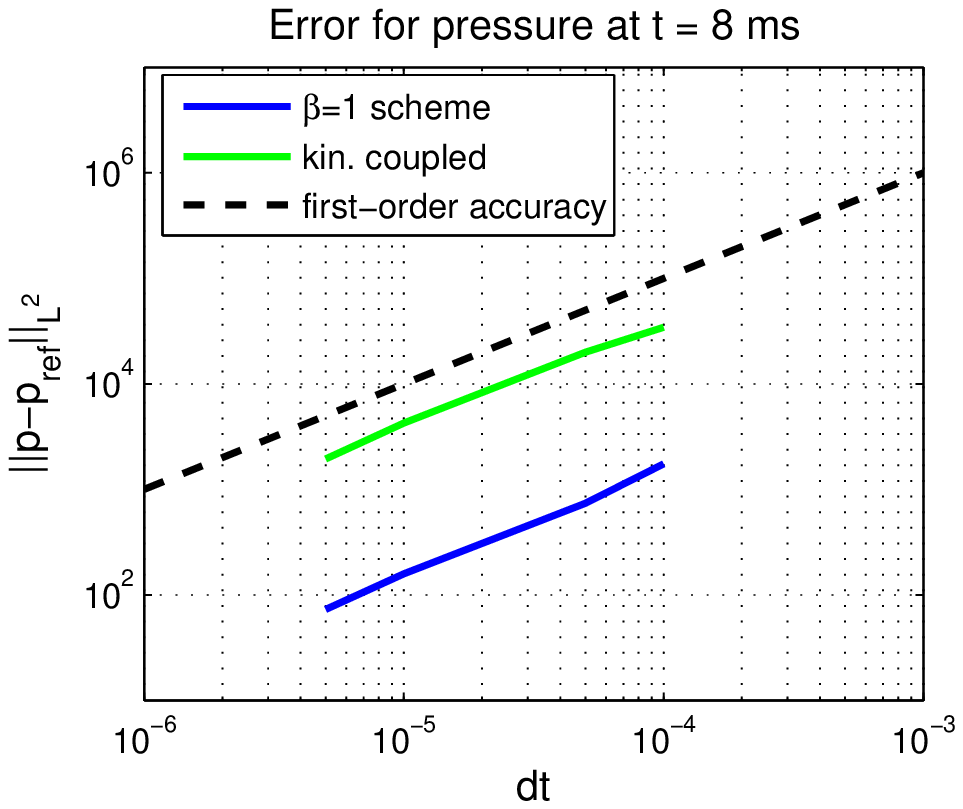} 
 \includegraphics[scale=0.65]{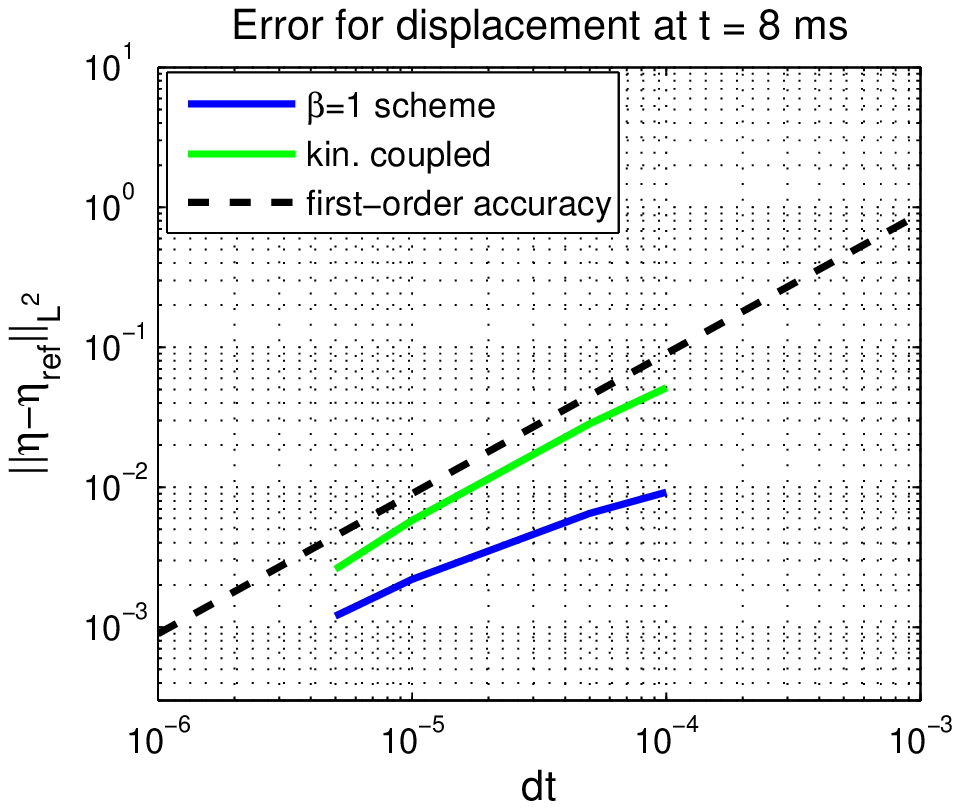}
 }
 \caption{Example 2: Figures show relative errors compared with the kinematically coupled scheme which is first-order accurate in time. Top left: Relative error for fluid velocity at t=8 ms. Top right: Relative error for fluid pressure at t=8 ms. Bottom: Relative error for displacement at t=8 ms.}
\label{error_ex2}
 \end{figure}

To study convergence in time we define the reference solution to be the one obtained with $\triangle t = 10 ^{-6}$, 
and we compute the $L^2$-norms of the difference in the  pressure, velocity and displacement 
between the reference solution and the solutions obtained with $\triangle t = 10^{-4}, 5 \times 10^{-5}, 10^{-5}$ and $5 \times 10^{-6}$.  
A comparison of the time convergence between our scheme and the kinematically coupled scheme is presented in Table~\ref{T7}. 
\begin{table}[htp!]
\begin{center}
{\scriptsize{
\begin{tabular}{| l  c  c  c  c  c  c |}
\hline
$ \triangle t $ & $||p-p_{ref}||_{L^2} $ & $L^2$ order & $||\boldsymbol u-\boldsymbol u_{ref}||_{L^2}$  &$L^2$ order & $ ||\boldsymbol \eta - \boldsymbol \eta_{ref}||_{L^2} $ & $L^2$ order \\
\hline
\hline
$ 10^{-4}$  & $1.75 \textrm{e}+03$  & - & $ 5.83  $& - & $0.0092$  & - \\
$  $  & $(3.422 \textrm{e}+04)$  & - & $ (91.175)  $& - & $(0.0513)$  & - \\
\hline
$5 \times 10^{-5}$ &  $739.9$  & 1.24 & $3.85$  & 0.6 & $  0.0065$  & 0.5\\
$ $ &  $(1.9943 \textrm{e}+04)$  & (0.78) & $(49.76)$  &(0.87) & $  (0.0284)$  & (0.85)\\
\hline
$ 10^{-5}$ & $ 158.46$  & 0.96 & $0.88$  & 0.92 & $ 0.0022$  & 0.66\\
$ $ &  $(4.24\textrm{e}+03)$  & (0.96) & $(10.016)$  &(0.99) & $  (0.0058)$  & (0.99)\\
\hline
$5 \times 10^{-6}$  & $72.67$  & 1.12 & $ 0.4$  &1.13 & $0.0012$  & 0.92 \\
$ $ &  $(1.94 \textrm{e}+03)$  & (1.12) & $(4.6)$  &(1.12) & $  (0.0026)$  & (1.16)\\
\hline
\end{tabular}
}}
\end{center}
\caption{Example 2: Convergence in time calculated at $t = 8$ ms. The numbers in the parenthesis show the convergence rate for the kinematically coupled scheme presented in~\cite{guidoboni2009stable}.}
\label{T7}
\end{table}
Figure~\ref{error_ex2} shows a comparison between the time convergence of our scheme and the kinematically coupled scheme. 
Similarly as before, our method implemented with $\beta = 1$ provides a gain in accuracy over the classical kinematically coupled scheme, which corresponds to $\beta = 0$.

\subsection{The Common Carotid Artery (CCA) Example}

\begin{figure}[ht!]
\centering{
 \includegraphics[scale=0.6]{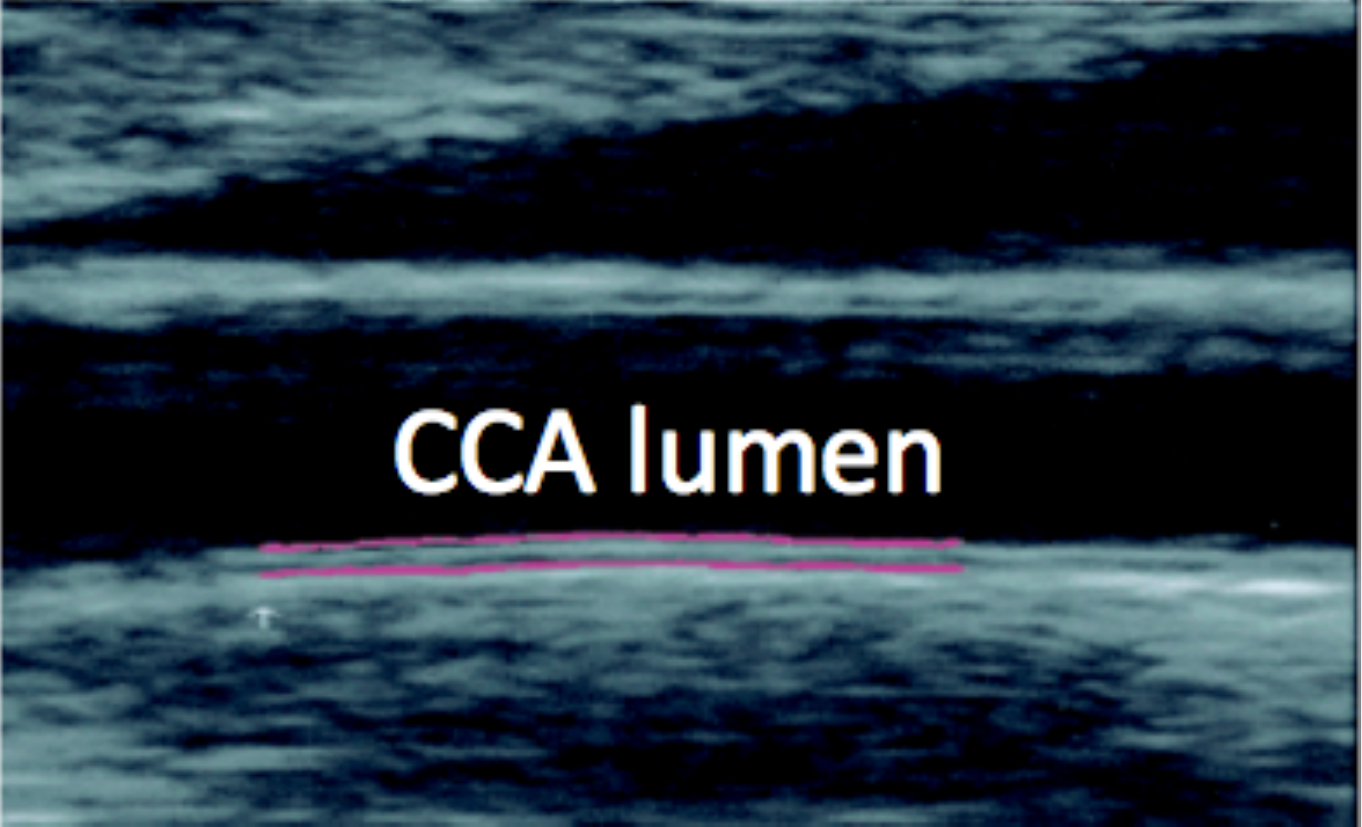}\hskip 0.5in
}
 \caption{A B-mode ultrasound image of CCA \cite{CCA_circulation} representing our computational domain. The black triangle above the CCA is the jugular vein (not a part of the computational domain). }
\label{CCA}
 \end{figure}

We conclude this manuscript by showing a simulation of flow and displacement
in the left common carotid artery.
A straight segment of the left CCA, such as the one shown in Figure~\ref{CCA}, is considered.
The geometric parameters, such as length, average radius, and wall thickness, are taken from the measurements reported
in~\cite{length,diameters, juonala2005risk,mokhtari2006differentiation,bussy2000intrinsic}, while the Youngs modulus 
is taken from the measurements reported in~\cite{mokhtari2006differentiation}.
Blood vessels are essentially incompressible and therefore have the Poisson's ratio of approximately 0.5~\cite{nichols2005mcdonald}. 
The Table in Figure~\ref{CCA} shows the values of the corresponding parameters that were used in our simulation.
They are within the ranges reported in the above-mentioned literature.
\begin{table}[ht]
\centering{
\begin{tabular}{|l l|}
\hline
\textbf{Parameters} & \textbf{CCA}  \\
\hline
\hline
\textbf{Radius} $H$ (cm)  & $0.3$  \\
\hline
%\textbf{Length} $L$ (cm) & $10$  \\
%\hline
\textbf{Fluid density} $\rho_f$ (g/cm$^3$)& $1.055$ \\
\hline
\textbf{Fluid viscosity} $\mu$ (g/(cm s)) & $0.04$    \\
\hline
\textbf{Wall density} $\rho_s $(g/cm$^3$) & $1.055$  \\
\hline
\textbf{Wall thickness} $h$ (cm) & $0.07$  \\
\hline
\textbf{Young's mod.} $E $(dynes/cm$^2$) & $2 \times 10^6$  \\
\hline
\textbf{Poisson's ratio} $\sigma $ & $0.5$   \\
\hline
\end{tabular}}
\caption{Geometry, fluid and structure parameters for the common carotid artery example.}
\label{table_CCA_param}
\end{table}

The structural viscosity constants $C_v$ and $D_v$  are equal to
$$C_v := 3 \times 10^4  \; \textrm{dynes/cm}^2 \cdot \textrm{s}, \quad D_v := C_v \; \sigma.$$
This choice of structural viscosity parameters was shown in~\cite{canic2006modeling} to resemble the viscous moduli of blood vessels. 
These parameters give rise to the values of the Koiter shell coefficients  given in Table~\ref{TCCAcoeff}.
\begin{table}[ht!]
\begin{center}
\begin{tabular}{| l  l  l  l  l|}
\hline
$C_0 = 1.7022 \times 10^6$ & $C_1=846.9$ & $C_2 = 3.1 \times 10^5 $ & $C_3=1.867 \times 10^5$  & $C_4=0$ \\
\hline
$D_0 = 23439.2$  & $D_1 = 9.527$  & $D_2 = 3500$ & $D_3 = 2100$ & $D_4=0$ \\
\hline
\end{tabular}
\end{center}
\caption{Koiter shell model coefficients for the common carotid artery example.}
\label{TCCAcoeff}
\end{table}

\begin{figure}[ht]
\centering{
\includegraphics[scale=1]{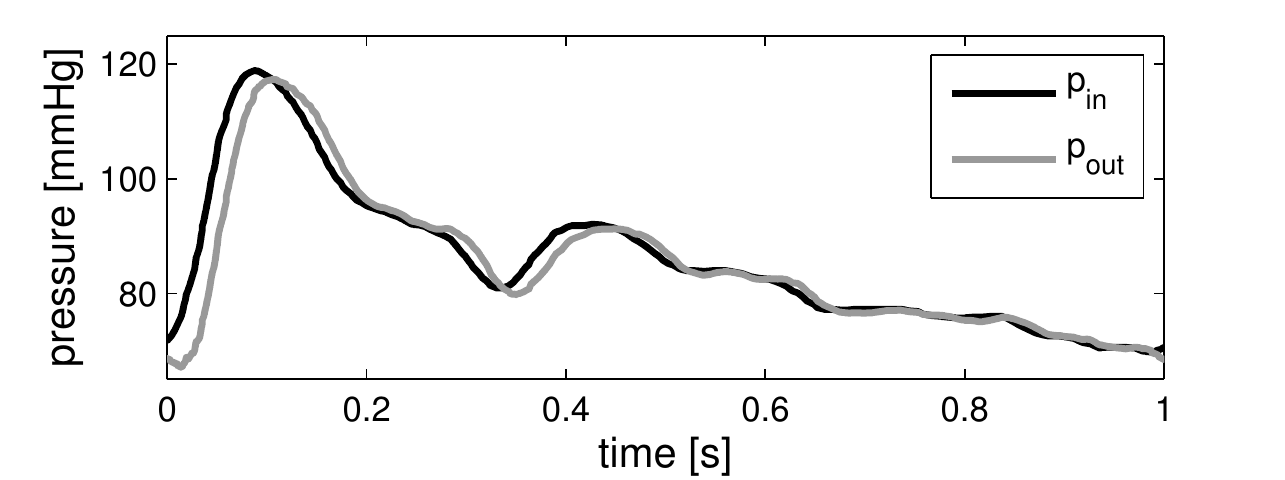}
}
\caption{The inlet and outlet pressure data. The average pressure drop is around 0.0673 mmHg per centimeter.}
\label{carotid_pressure_data}
\end{figure}

\begin{figure}[ht!]
\centering{
 \includegraphics[scale=0.95]{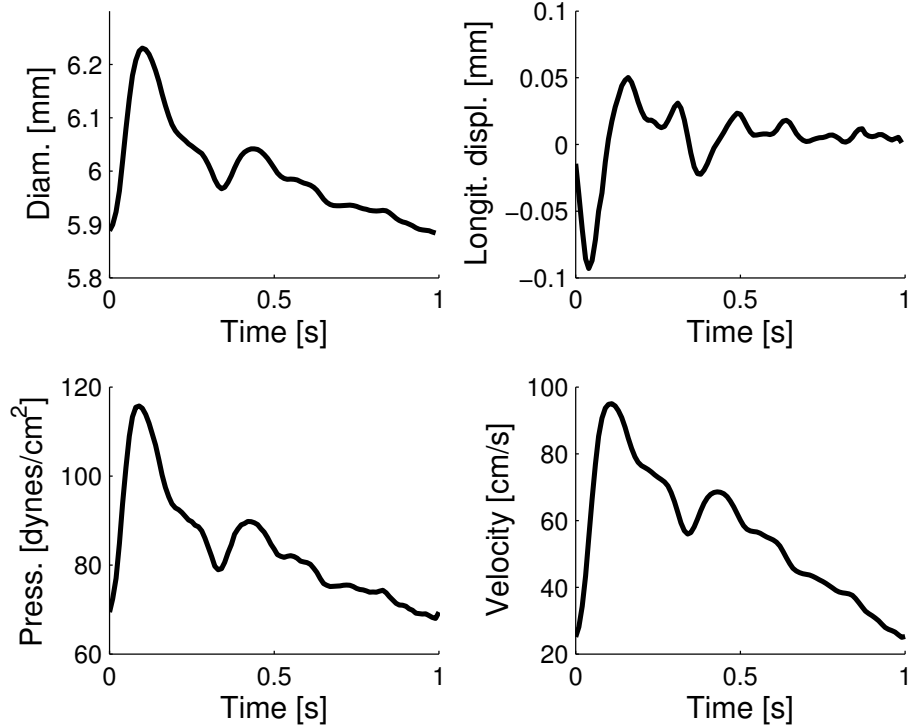}\hskip 0.5in
 }
 \caption{From top left, to bottom right: radial displacement, longitudinal displacement, pressure, and velocity, calculated at
 the mid-point of the CCA segment, vs. time.}
\label{CCA_time}
 \end{figure}

We study blood flow driven by the inlet and outlet pressure data shown in Figure~\ref{carotid_pressure_data}.
The morphology of the pressure wave in the left CCA was obtained from \cite{warriner2008viscoelastic}.
The average pressure drop equals around 0.0673 mmHg per centimeter, which produces the local Reynolds number of around 1000.
This is associated with the maximum blood flow velocity of around 100 cm/s, which is typical 
for CCA~\cite{rohren2003spectrum, bouthier1985pulsed, dizaji2009estimation, lee1999assessment}.
Indeed, the results of our numerical simulation, presented in Figure~\ref{CCA_time}, show the velocity 
ranging between 22 cm/s and 97 cm/s, which is within the expected values for the CCA. 

Radial wall displacement has been well examined by many experimental studies~\cite{bussy2000intrinsic, dammers2003shear, meinders2003assessment, samijo1998wall, dizaji2009estimation}. Maximum radial displacement decreases with age, and usually varies between 0.1 mm and 0.38 mm,
 i.e. between 3\%and 13\%  of the vessel's radius. 
 Indeed, our simulation, shown in Figure~\ref{CCA_time} top left, indicates maximum radial displacement around 6\%, which is 
well within the normal range. 

Figure~\ref{CCA_time} top right,  shows longitudinal displacement computed using our thin shell model. 
We see that the longitudinal displacement in our simulations lies between
${\rm min} \ \eta_z=  -0.1$ mm and ${\rm max} \ \eta_z = 0.05$ mm, which implies the total longitudinal displacement 
(defined to be  ${\rm max} \ \eta_z - {\rm min} \ \eta_z$) of 0.15mm.
This is in good agreement with experimental studies obtained using B-mode ultrasound speckle tracking method
and/or B-mode ultrasound velocity vector imaging~\cite{Svedlund,cinthio2006longitudinal, persson2003new,svedlund2011longitudinal},
which report the total longitudinal displacement in a healthy CCA between 0.052 mm and 0.302 mm.

 \begin{figure}[ht]
\centering{
\includegraphics[scale=1.0]{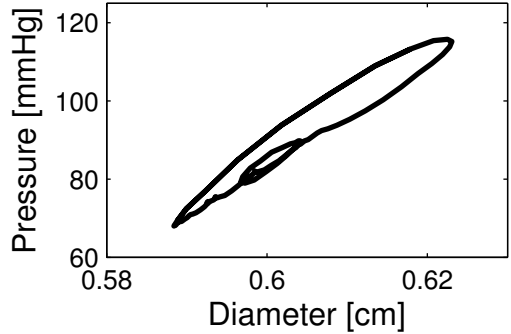}
}
\caption{Hysteresis between vessel wall diameter and pressure at the center of the vessel, over a cardiac cycle.}
\label{F532}
\end{figure}

Finally, we observe the captured viscoelastic properties. 
The viscoelastic effect is  visible in the stress-strain relationship of the arterial wall, which exhibits hysteresis. Figure~\ref{F532} shows hysteresis between the vessel diameter and pressure at the center of the vessel over one cardiac cycle.
Hysteresis can be quantified by the energy dissipation ratio (EDR), which is a measure of the area inside the diameter-pressure loop relative to the measure of areas inside and under the loop, i.e.,  $EDR = A_1/(A_1 + A_2) \times 100$\%.
Walls with higher viscoelasticity have larger area inside the loop, resulting in higher EDR. In our simulations EDR is $8.5 \%$. This is comparable to the results in~\cite{warriner2008viscoelastic} which show EDR of $7.8 \%$ for a young subject.

 \section{Conclusions}
In this manuscript we proposed a new thin structure model capturing radial and longitudinal displacement of arterial walls,
and have designed a modification of a loosely coupled partitioned scheme (the kinematically coupled scheme \cite{guidoboni2009stable}) to
numerically simulate the resulting fluid-structure interaction problem between blood flow and arterial walls.
The proposed arterial wall model is given by the linearly viscoelastic, cylindrical Koiter shell model. The fluid and structure are fully
coupled using the kinematic and dynamic coupling conditions. The new loosely coupled scheme
(the kinematically coupled $\beta$-scheme),
which is 1st-order accurate in time, and 2nd-order accurate in space, is based on a modified
Lie splitting. In \cite{stability} it is shown that this scheme is unconditionally stable for all $0\le \beta \le 1$.
Two test problems were presented showing that the kinematically coupled $\beta$-scheme has accuracy which is 
comparable to that of the monolithic scheme by Badia, Quaini, and Quarteroni \cite{quaini2009algorithms,badia2008splitting}
while retaining the main advantages of loosely coupled partitioned schemes such as
modularity, easy implementation, and low computational costs (no sub-iterations between the 
fluid and structure sub-solvers are necessary for convergence). We believe that the main reason for the increase 
in accuracy is related to the stronger coupling between the leading effect of the fluid load, given by the fluid pressure,
and structure elastodynamics, which are now strongly coupled in the last step (Step 3) of the splitting.
It remains to be investigated how does the choice of $\beta$, which would provide highest possible accuracy
of the kinematically coupled $\beta$ scheme, depend on the coefficients in the problem. 

Although the methodology discussed in this manuscript was presented for 2D problems, 
there is nothing is the proposed time-splitting scheme that depends on the dimension of the spatial domain.
The same methodology can be applied to 3D problems, while, of course, the implementation of the proposed
algorithm in 3D would be much more complex in that case.

\if 1 = 0
Application of the proposed model and scheme to a set of physiologically relevant problems is currently investigated.
In the work near completion \cite{physiological} we use the proposed numerical scheme to recreate the
radial and longitudional displacement in the common carotid artery that was measured using
{\sl in vivo} ultrasound speckle tracking techniques \cite{warriner2008viscoelastic,LyonLongitudinalMovie}. Our preliminary results show excellent agreement
with experimental measurements not only in the radial and longitudinal displacement, but also in
the hysteresis behavior due to viscoelasticity of vessel walls.
\fi

Extensions of the proposed methodology to study FSI between an incompressible, viscous fluid and a {\bf thick} structure (arterial walls), modeled using
equations of 2D/3D elasticity (e.g., St. Venant-Kirchhoff model, used, e.g., in \cite{Bazilevs1,Bazilevs2} to model arterial walls), are under way. 

The research presented in this manuscript provides a first step in our effort to capture multi-layered structure of arterial walls
and their interaction with blood flow. In modeling the intima-media/adventitia complex, the coupling between
a thin shell (intima) allowing radial and longitudinal displacement, and a thick structure (media/adventitia)
is important. 
Development of the model presented in this manuscript is a crucial first step.
Our preliminary results show that the modified kinematically coupled scheme 
proposed in this manuscript
is perfect for the numerical solution of such a complex, multi-physics FSI problem.
Research in this direction in under way.

 \section{Acknowledgements}
The authors would like to thank Boris Muha for his help with the manuscript. 
Buka\v{c} acknowledges partial graduate student support by the National Science Foundation through grants DMS-1109189 and DMS-0806941.
\v{C}ani\'c acknowledges partial research support by NSF under grants DMS-1109189, DMS-0806941, by the Texas Higher Education Board under
Advanced Research Program,  Mathematics-003652-0023-2009, and by the joint NSF and  NIH support under grant DMS-0443826.
Glowinski acknowledges partial research support by NSF under grants DMS-0811153, DMS-0914788 and DMS-0913982. 
Quaini acknowledges partial research support by NSF under grant DMS-1109189, and partial post-doctoral research support
by the Texas Higher Education Board under grant ARP Mathematics-003652-0023-2009.

\bibliographystyle{model1-num-names}
\bibliography{method}
\end{document}